\documentclass[11pt]{article}
\usepackage[margin=1in]{geometry}
\usepackage{amsmath,amssymb,amsthm,hyperref,microtype}
\hypersetup{colorlinks=true,linkcolor=blue,citecolor=blue,urlcolor=blue}
\newtheorem{theorem}{Theorem}[section]
\newtheorem{lemma}[theorem]{Lemma}
\newtheorem{corollary}[theorem]{Corollary}
\newtheorem{proposition}[theorem]{Proposition}

\theoremstyle{definition}
\theoremstyle{remark}\newtheorem{remark}[theorem]{Remark}
\title{Induced Independence Levels and Signed Enumeration in Powers of Cycles}
\author{Jing Zhang\\Department of Mathematics\\Governors State University\\University Park, Illinois 60484, USA\\\texttt{jzhang@govst.edu}\and Yuan Li\thanks{Corresponding author. E-mail: \texttt{liyu@wssu.edu}.}\\Department of Mathematics\\Winston-Salem State University\\Winston-Salem, North Carolina 27110, USA\\\texttt{liyu@wssu.edu}}
\date{}
\begin{document}\maketitle
\begin{abstract}
For a graph $G$ we study the distribution of the independence numbers
$\alpha(G[S])$ over its vertex-induced subgraphs, together with signed
Boolean-lattice sums over fixed independence-number levels.
For the $r$th power $C_n^r$ of a cycle, we establish, in the range
$n\ge 2rk+r+1$, a unique representation of each subset $S$ with
$\alpha(C_n^r[S])=k$ by an independent $k$-set and its forward
neighborhood. The resulting Boolean-interval partition gives exact
enumeration by subset cardinality. We also prove reconstruction of
partially deleted forward neighborhoods and determine the support,
signs, and multiplicities of the corresponding upper M\"obius transform.
For cycle squares we develop a complementary greedy-chain method for
signed enumeration, with a five-transition classification and explicit
staircase and two-gap cancellation identities.
\end{abstract}

\section{Introduction}
\subsection{Induced independence levels}
Let $G=(V,E)$ be a finite simple graph and let $G[S]$ denote the
subgraph induced by $S\subseteq V$. The independence number
$\alpha(G[S])$ measures the largest set of mutually nonadjacent
vertices within $S$. A natural induced-subgraph enumerator is
\[
 F_G(y,t)=\sum_{S\subseteq V}y^{|S|}t^{\alpha(G[S])}.
\]
Its coefficients
\[
 A_{j,k}(G)=\#\{S\subseteq V:|S|=j,\ \alpha(G[S])=k\}
\]
record the distribution of the induced independence number at
every subset cardinality. The familiar independence polynomial
$I_G(z)=\sum_k i_k(G)z^k$ counts only those $S$ that are themselves
independent. By contrast, $A_{j,k}(G)$ includes all induced
subgraphs, including those with many edges.

We focus on powers of cycles. On $\mathbb Z_n$, the graph $C_n^r$
joins vertices whose cyclic distance is at most $r$. We take the
usual nondegenerate range of parameters in the results below.
The case $r=2$ admits additional descriptions by interlacing
greedy chains. For signed enumeration we study alternating sums over subsets
of a fixed independence level and their Boolean-lattice
M\"obius transforms. Such sums can exhibit cancellations that
are invisible in the nonnegative coefficients of $F_G$.
Our methods exploit two complementary descriptions of the
independence levels: a uniform forward-neighborhood
decomposition in the stable range for $C_n^r$, and local
greedy-chain and canonical-requirement decompositions for
$C_n^2$.

\subsection{Related literature}
\paragraph{Independent sets and powers of cycles.}
The independence polynomial and its combinatorial interpretations
are standard objects; Levit and Mandrescu~\cite{LevitMandrescu2005}
provide a survey. Brown and Hoshino~\cite{BrownHoshino2009}
obtain independence polynomials for a family of circulants
including powers of cycles. In particular, the number of
independent $k$-sets of $C_n^r$ is the classical expression
\[
 i_k(C_n^r)=\frac{n}{n-rk}\binom{n-rk}{k}
\]
when the parameters permit such sets. The underlying circular
spacing enumeration goes back to classical work of
Kaplansky~\cite{Kaplansky1943} and Yamamoto~\cite{Yamamoto1956};
the cycle-power formulation is also treated by Brown and Hoshino.
The formula is used in our enumeration but is not claimed as a new result.
Codara and D'Antona~\cite{CodaraDAntona2012} also examine
independent subsets of powers of paths and cycles. Their
work concerns the classical independent-set enumeration,
whereas the level-set partition below describes all induced
subsets with a prescribed independence number.
White~\cite{White2026} studies cyclic sieving for independent
sets of graphs, including cycle powers, and places these
counts in a symmetry-sensitive setting.
The distinctive question here is not how many independent
$k$-sets exist, but how \emph{all} vertex subsets with induced
independence number $k$ can be represented and counted.

Bermond and Peyrat~\cite{BermondPeyrat1989} study induced
subgraphs of powers of cycles under minimum-degree and
density constraints. Their work provides direct context
for structural questions about induced subgraphs of this
graph family, though their principal constraints differ
from the independence-level statistic considered here.

\paragraph{Induced-subgraph generating functions.}
Tittmann, Averbouch, and Makowsky~\cite{TittmannAverbouchMakowsky2011}
study a bivariate graph polynomial enumerating induced
subgraphs by order and number of connected components.
Their approach illustrates the wider program of retaining
structural information across all induced subgraphs rather
than only counting a distinguished class of subsets.
The statistic in $F_G$ is instead $\alpha(G[S])$; the
component-counting polynomial does not directly supply
our formulas.

A stronger informational comparison comes from the
Kromatic symmetric function. Pierson~\cite{Pierson2024}
uses this invariant to count specified classes of induced
subgraphs, including small induced configurations and
stars together with isolated vertices. In subsequent work,
Pierson~\cite{Pierson2025} proves that the Kromatic
symmetric function determines the multiset of independence
polynomials of all vertex-induced subgraphs; see also
Pierson and Samanta~\cite{PiersonSamanta2026}.
For nonempty $S$, the degree of $I_{G[S]}(z)$ is
$\alpha(G[S])$ and its coefficient of $z$ is $|S|$.
The empty set is recovered separately. Hence this multiset
determines $F_G(y,t)$. Thus our goal is not to introduce
a graph invariant carrying information unavailable from
existing symmetric-function invariants. Instead, we seek
explicit structural and signed enumeration formulas for
powers of cycles that are not furnished by this general
recoverability statement.

\paragraph{Boolean intervals and signed enumeration.}
Inclusion--exclusion and M\"obius inversion on the Boolean
lattice are classical; see Stanley~\cite{StanleyEC1}.
For a family of subsets specified by an induced
independence number, the challenge is to make the
alternating sums explicit. A partition into Boolean
intervals provides one mechanism, while canonical
requirements for independent sets and greedy-chain
restrictions provide another. The general inversion
principle is standard; the interval structure and the
cycle-square cancellation identities established here
are the objects of study.

\subsection{Contributions and organization}
Our first principal contribution is a \emph{unique Boolean-interval
partition of an entire induced independence-number level},
not merely the classical count of independent sets.
Specifically, we establish a forward-slots
representation for $C_n^r$ in the stable range
$n\ge2rk+r+1$. For an independent $k$-set $I$ put
\[
 Q_r(I)=\bigcup_{x\in I}\{x,x+1,\ldots,x+r\}
 \quad\text{(indices modulo $n$)}.
\]
We prove that every $S$ with $\alpha(C_n^r[S])=k$
has a unique representation
$I\subseteq S\subseteq Q_r(I)$.
In particular, the intervals $[I,Q_r(I)]$ partition
the entire independence level and give
\[
 A_{j,k}(C_n^r)
 =\frac{n}{n-rk}\binom{n-rk}{k}
   \binom{rk}{j-k}.
\]
The independent-set factor in this expression is known;
the interval partition is the structural statement that
explains the full induced-level enumeration. We further
prove injectivity for reconstruction of a forward-neighborhood
set after deleting selected slot starts. As a consequence
of the partition and this reconstruction theorem, we
determine the exact support, signs, and multiplicities
of the upper Boolean-lattice M\"obius transform of
the level indicator. The M\"obius formula is presented
as a consequence of the structural theorems, not as a
new general inversion principle.

Our second principal contribution concerns the more
delicate signed enumeration for $C_n^2$ outside the
scope of a single uniform level count. We use
interlacing extremal greedy chains to constrain
possible local transitions and to identify
canonical requirements for maximum independent sets.
This yields signed weights, staircase identities,
boundary cancellation, and classifications of
interior two-gap configurations. The greedy-chain
and forward-slots approaches share the same
underlying induced-independence statistic, but
address different structural and signed questions.

For a guide to the logical structure, the general structural result is
Theorem~\ref{thm:slots-bijection}; the reconstruction
Lemma~\ref{lem:partial-slot-injectivity} and the M\"obius
Theorem~\ref{thm:stable-mobius-distribution} build on that
partition. In the separate cycle-square theory,
Theorem~\ref{thm:general-canonical-requirements} supplies the
uniform five-transition signed-weight mechanism. The subsequent
staircase, boundary, and two-gap results identify specific
configurations and coefficients; their additional classifications
are not claimed to follow from the uniform weight formula alone.

The paper first develops the cycle-square greedy-chain
method and its applications. It then presents the
uniform cycle-power forward-slots theory and its
M\"obius consequences. The final sections describe
restricted transforms and the associated
simplicial-complex interpretation. These standard
algebraic-topological tools are included to clarify
the signed combinatorics, rather than as independent
claims of novelty.

\section{Greedy-chain signed enumeration and staircase theory}
\noindent\textbf{Why the boundary scores matter.}
The cycle square differs from the path square by three wraparound
edges. We first retain the edge $\{1,n\}$, writing
$G_n=P_n^2+\{\{1,n\}\}$, and distinguish two boundary patterns
that are independent in $G_n$ but cease to be independent when
$\{1,n-1\}$ and $\{2,n\}$ are inserted. Their maximum
attainable cardinalities are denoted by $U(S)$ and $V(S)$;
$X(S)$ is the largest cardinality among all other boundary patterns
(the precise definitions follow below). The condition
$U(S)=V(S)=k>X(S)$ therefore isolates subsets for which
both forbidden patterns attain the maximum, whereas every
cyclically admissible pattern is smaller. We call these
\emph{both-active deficits}. Their signed contribution is
recorded by the polynomial $B_n^{\mathrm{both}}(t)$ defined below.
The following results evaluate specific coefficients and families
of these contributions; no unstated global identity is needed.

\subsection{An exact signed forced-endpoint identity}
\label{subsec:greedy-boundary}
We begin with a forced-endpoint signed sum for path squares.
This result provides a direct example of how greedy chains and
Boolean intervals yield cancellation, and it explains the
four-step threshold that reappears in the boundary analysis.

Let $P_n^2$ be the square of the path on $[n]$.  Say that a vertex of an
induced subgraph is \emph{forced} if it belongs to every maximum independent
set.  Put
\[
 \mathcal B_n(t)=\sum_{\substack{S\subseteq[n]\\
  1,n\text{ forced in }P_n^2[S]}}
 (-1)^{n-|S|}t^{\alpha(P_n^2[S])-1}.
\]
\begin{theorem}[Exact signed forced-endpoint distribution]
\label{thm:greedy-boundary}
For every $n\ge7$,
\[
 \mathcal B_n(t)=(-1)^n\sum_{j=1}^{\lfloor(n-2)/4\rfloor}t^j.
\]
More precisely, let $N_k(n)$ denote the signed count of admissible pairs
of left- and right-greedy chains of length $k$, as defined in the proof.
Then for $k\ge2$,
\[
 N_k(n)=\begin{cases}1,&n\ge4k-2,\\0,&n<4k-2.\end{cases}
\]
\end{theorem}
\begin{proof}
For an available set $S\subseteq[n]$, left greedy repeatedly chooses the
smallest available vertex at least three positions after its preceding
choice.  Right greedy is its reversed counterpart.  Both give maximum
independent sets: replacing the first element of any independent set by
the smallest available element does not decrease its size, and induction
finishes the argument.  Vertex $1$ is forced exactly when right greedy
chooses it; similarly $n$ is forced exactly when left greedy chooses it.

For a set contributing to $\mathcal B_n$, write the two greedy chains in
increasing order as
\[
 L=\{a_1<\cdots<a_k\},\qquad R=\{b_1<\cdots<b_k\},
 \qquad a_1=b_1=1,\quad a_k=b_k=n.
\]
All successive gaps within either chain are at least three.  Define
\[
 A(L)=[n]\cap\bigcup_{j=1}^k[a_j,a_j+2],\qquad
 B(R)=[n]\cap\bigcup_{j=1}^k[b_j-2,b_j].
\]
The sets having these two greedy chains form exactly the Boolean interval
\[
 L\cup R\ \subseteq\ S\ \subseteq\ A(L)\cap B(R).
\]
Indeed, each available vertex must lie in the two indicated greedy
windows; conversely, the windows prevent an earlier eligible successor
than the specified next greedy choice.  Every set in this interval has
independence number $k$.  Summing $(-1)^{|S|}$ over the interval cancels
unless $L\cup R=A(L)\cap B(R)$, in which case the signed contribution
is $(-1)^{|L\cup R|}$.  Call such a pair \emph{admissible} and write
\[
 N_k(n)=\sum_{\substack{(L,R)\text{ admissible}\\|L|=|R|=k}}
 (-1)^{|L\cup R|}.
\]
It follows that $\mathcal B_n(t)=(-1)^n\sum_{k\ge2}N_k(n)t^{k-1}$.

We characterize admissibility.  Set
$A_j=a_{j+1}-a_j$ and $\delta_j=b_j-a_j$.
The inclusion $L\cup R\subseteq A(L)\cap B(R)$ forces each $a_j$ to
lie within two positions to the left of some $b_i$, and each $b_i$ to
lie within two positions to the right of some $a_j$.  Because the
windows for each chain are disjoint and preserve order, the resulting
matching of chain vertices is the identity.  Thus
$0\le\delta_j\le2$.
If $\delta_j=2$, the point $a_j+1$ belongs to both window unions but
not to $L\cup R$, contradicting admissibility.  Hence
$\delta_j\in\{0,1\}$.

For these offsets, the diagonal intersection of the two windows at
index $j$ is $[a_j,b_j]$, consisting only of chain vertices.  Windows
with indices differing by at least two do not intersect.  The remaining
possible intersection is
$[a_j,a_j+2]\cap[b_{j+1}-2,b_{j+1}]$.
It is nonempty exactly when
$A_j+\delta_{j+1}\le4$; in that case it contains $a_j+2$, which
belongs to neither chain.  We conclude that admissibility is equivalent
to
\begin{equation}\label{eq:greedy-gap-criterion}
 \delta_j\in\{0,1\},\qquad A_j+\delta_{j+1}\ge5
 \quad(1\le j<k).
\end{equation}
The endpoint offsets satisfy $\delta_1=\delta_k=0$.

Put $r=k-2$, and let $m$ of the $r$ interior offsets be one.  Then
$|L\cup R|=k+m$, so the sign is $(-1)^{k+m}=(-1)^{r-m}$.
Write $A_j=5-\delta_{j+1}+f_j$ with $f_j\ge0$.
Since $\sum_j A_j=n-1$, the slack variables satisfy
\[
 \sum_{j=1}^{k-1}f_j=n-5k+4+m.
\]
For fixed offsets, the number of choices is
$\binom{n-4k+2+m}{k-2}$ whenever the slack sum is nonnegative,
and zero otherwise.  If $n<4k-2$, the slack sum is negative for every
$m\le k-2$, hence $N_k(n)=0$.

Suppose $n\ge4k-2$, and set $N=n-4k+2\ge0$.
With the usual convention $\binom qr=0$ for $0\le q<r$, the signed
count is
\[
\begin{aligned}
 N_k(n)
 &=\sum_{m=0}^{r}(-1)^{r-m}\binom rm\binom{N+m}{r}\\
 &=[x^r](1+x)^N\sum_{m=0}^{r}(-1)^{r-m}\binom rm(1+x)^m\\
 &=[x^r](1+x)^Nx^r=1.
\end{aligned}
\]
The asserted formula for $\mathcal B_n$ follows upon writing $j=k-1$.
\end{proof}

\begin{remark}[The four-step threshold]
The gap criterion \eqref{eq:greedy-gap-criterion} requires at least four
positions per interior gap and five for the final gap, since
$\delta_k=0$.  Hence
$n-1=\sum_{j=1}^{k-1}A_j\ge4(k-2)+5=4k-3$.
The threshold $n\ge4k-2$ is therefore geometric, while the signed
coefficient $1$ is explained by a finite-difference cancellation.
The threshold is a consequence of the gap geometry, while the
signed coefficient follows from a finite-difference identity.
\end{remark}

\begin{lemma}[Boundary forcing in the both-active case]
\label{lem:both-active-boundary-forcing}
Let $n\ge9$, $S\subseteq[n]$, and $G_n=P_n^2+\{\{1,n\}\}$.
For $i\in\{0,1,2\}$ and $j\in\{0,n-1,n\}$, let
$q_{ij}(S)$ be the maximum cardinality of an independent set in
$G_n[S]$ whose intersections with $\{1,2\}$ and $\{n-1,n\}$
are respectively $\varnothing$ or $\{i\}$, and $\varnothing$ or
$\{j\}$ (with $0$ denoting the empty choice). Set $q_{ij}=-\infty$
when the choice is impossible, and put
\[
 U=q_{1,n-1},\qquad V=q_{2,n},\qquad
 X=\max\{q_{ij}:(i,j)\ne(1,n-1),(2,n)\}.
\]
If $U=V=k>X$, then
\[
 \{1,2,4,n-3,n-1,n\}\subseteq S.
\]
Writing $a(p,q)=\alpha(P_n^2[S\cap[p,q]])$, we also have
\[
 a(4,n-4)=a(5,n-3)=k-2,\qquad a(5,n-4)=k-3.
\]
In particular, $4$ belongs to every maximum independent set of
$P_n^2[S\cap[4,n-4]]$, and $n-3$ belongs to every maximum independent
set of $P_n^2[S\cap[5,n-3]]$.
\end{lemma}
\begin{proof}
Finiteness of $U$ and $V$ requires $1,2,n-1,n\in S$. The
corresponding boundary choices exclude, respectively, vertices
$3,n-3,n-2$ and $3,4,n-2$ from the interior. Consequently
\[
 U=2+a(4,n-4),\qquad V=2+a(5,n-3).
\]
The choice $(2,n-1)$ is compatible with $G_n$, so
\[
 X\ge q_{2,n-1}=2+a(5,n-4).
\]
Since $U=V=k>X$, these relations imply
\[
 a(4,n-4)=a(5,n-3)=k-2,\qquad a(5,n-4)\le k-3.
\]
Deleting vertex $4$ from the first interval, or vertex $n-3$
from the second, decreases the independence number by at most one.
Thus $a(5,n-4)\ge k-3$, proving equality. Both deletions therefore
lower the independence number by exactly one. Each deleted vertex
must be present and belong to every maximum independent set of its
respective interval. This proves the assertions.
\end{proof}

\begin{lemma}[Boundary scores and independence numbers]
\label{lem:both-active-identifications}
Let $U(S)=V(S)=k>X(S)$, with the boundary scores of
Lemma~\ref{lem:both-active-boundary-forcing}. Then
$\alpha(C_n^2[S])=X(S)=k-1$ and
$\alpha(P_n^2[S])=\alpha(G_n[S])=k$.
\end{lemma}
\begin{proof}
$C_n^2$ is $P_n^2$ plus the three wraparound edges
$\{1,n\}$, $\{1,n-1\}$, $\{2,n\}$; $G_n$ is $P_n^2$ plus
$\{1,n\}$ only. A $C_n^2$-independent set is
$G_n$-independent with a non-forbidden boundary pattern,
so $\alpha(C_n^2[S])\le X(S)$; conversely a
$G_n$-independent set with a non-forbidden pattern contains
neither $\{1,n-1\}$ nor $\{2,n\}$, the only $C_n^2$-edges
missing from $G_n$, so it is $C_n^2$-independent and
$\alpha(C_n^2[S])=X(S)$. Deleting $1$ from a
$U$-realising $k$-set gives a $C_n^2$-independent
$(k-1)$-set, so $X\ge k-1$, and $X<k$ by hypothesis; hence
$X=k-1$. Every $G_n$-independent set has some boundary
pattern, so its size is at most $\max(U,V,X)=k$, and $U$
attains $k$, so $\alpha(G_n[S])=k$. Finally
$\alpha(P_n^2[S])\ge k$ by inclusion, and a
$P_n^2$-independent set violates at most one wraparound
edge, so deleting one endpoint leaves a $C_n^2$-independent
set of size at most $k-1$; hence $\alpha(P_n^2[S])\le k$.
\end{proof}

\begin{corollary}[Greedy-chain endpoint rigidity in the both-active case]
\label{cor:both-active-endpoint-rigidity}
Under the hypotheses of Lemma~\ref{lem:both-active-boundary-forcing},
let $L=(a_1<\cdots<a_k)$ and $R=(b_1<\cdots<b_k)$ be the
left- and right-greedy maximum independent sets of $P_n^2[S]$.
Then
\[
 a_1=1,\qquad a_2=4,\qquad b_1=2,\qquad
 b_{k-1}=n-3,\qquad a_k=n-1,\qquad b_k=n.
\]
In particular, the first gap of $L$ and the last gap of $R$ both
equal three.
\end{corollary}
\begin{proof}
Since $\alpha(G_n[S])=k$ and, by Lemma~\ref{lem:both-active-identifications}, $\alpha(C_n^2[S])=X=k-1$,
the independence number of $P_n^2[S]$ is $k$: it is at least $k$
by graph inclusion, and at most $k$ because a maximum independent
set of $P_n^2[S]$ can violate at most one of the three wraparound
edges $\{1,n\}$, $\{1,n-1\}$, $\{2,n\}$. Indeed, any two of
these edges have endpoints that are adjacent in $P_n^2$, so an
independent set cannot contain both edges; deleting one endpoint
of the unique violated edge leaves an independent set of
$C_n^2[S]$.

The preceding lemma gives $1,2,4,n-3,n-1,n\in S$.
The left-greedy chain consequently begins with $a_1=1$, $a_2=4$,
and the right-greedy chain ends with $b_{k-1}=n-3$, $b_k=n$.
Because $V=k$, there is an independent $k$-set containing $2$ and
$n$. For independent $k$-sets in a path square, the right-greedy
chain is coordinatewise at least as large as any other such set:
indeed $b_k=n$ is maximal, and $b_j$ is the largest element of
$S$ that is at most $b_{j+1}-3$, whereas the $j$-th element $c_j$
of any other maximum independent set satisfies
$c_j\le c_{j+1}-3\le b_{j+1}-3$ by descending induction;
therefore $b_1\ge2$. On the other hand, since $a_1=1$, the
right-greedy chain's first element is at most $3$: if $b_1\ge4$,
then $1,b_1,\ldots,b_k$ would be an independent $(k+1)$-set.
Thus $b_1\in\{2,3\}$.

If $b_1=3$, the independent $(k-1)$-set
$\{b_1,\ldots,b_{k-1}\}$ lies in $S\cap[3,n-3]$.
Together with $n$, it gives an allowed boundary choice $(0,n)$
of size $k$ in $G_n[S]$, contradicting $X<k$.
Hence $b_1=2$.
Reflecting the path by $i\mapsto n+1-i$ exchanges the two
forbidden boundary choices and the left- and right-greedy chains.
The proved assertion $b_1=2$ for the reflected configuration
therefore gives $a_k=n-1$.
\end{proof}

\begin{proposition}[Cancellation at the maximal both-active chain length]
\label{prop:both-active-upper-cancellation}
Let $n\ge9$ and let $S\subseteq[n]$ satisfy $U(S)=V(S)=k>X(S)$,
with the boundary scores as in
Lemma~\ref{lem:both-active-boundary-forcing}. Then $n\ge3k-1$.
For each fixed $k\ge 4$, if $n=3k-1$, the signed sum over all
such subsets is zero:
\[
 \sum_{\substack{S\subseteq[n]\\U(S)=V(S)=k>X(S)}}
 (-1)^{n-|S|}t^{\alpha(C_n^2[S])}=0.
\]
\end{proposition}
\begin{proof}
By Corollary~\ref{cor:both-active-endpoint-rigidity}, the
left-greedy chain $L=(a_1<\cdots<a_k)$ has $a_1=1$ and $a_k=n-1$.
Its successive gaps are at least three, so
\[
 n-2=a_k-a_1\ge3(k-1),
\]
which proves $n\ge3k-1$.

Suppose equality holds. All gaps in $L$ equal three. The
right-greedy chain $R=(b_1<\cdots<b_k)$ likewise has $b_1=2$,
$b_k=n$, and all its gaps equal three. Hence the chains are uniquely
specified by
\[
 L=\{1,4,7,\ldots,n-1\},\qquad
 R=\{2,5,8,\ldots,n\}.
\]
Put $S_0=L\cup R$ and
\[
 \mathcal J=\{S:\ S_0\subseteq S\subseteq A(L)\cap B(R)\},
\]
where $A(L)$ and $B(R)$ are the greedy coverage sets defined above.
The Boolean-interval lemma shows that $\mathcal J$ is precisely the
family of subsets having these two greedy chains. Moreover, for every
$S\in\mathcal J$, the path-square independence number equals $k$.
Although $L$ is independent in the path square, its cyclic endpoints
$1,n-1$ are adjacent in the cycle square. Since
$\alpha(C_n^2)=\lfloor n/3\rfloor=k-1$, we have
$\alpha(C_n^2[S])\le k-1$. Deleting vertex $1$ from $L$ gives an
independent $(k-1)$-set of $C_n^2[S]$, so equality holds. Thus every
$S\in\mathcal J$ has deficit one.

The vertex $3$ belongs to $A(L)$ (it is covered by $1$) and to
$B(R)$ (it is covered by $5$), but not to $S_0$. Consequently
$S\mapsto S\mathbin{\triangle}\{3\}$ is a fixed-point-free involution
of $\mathcal J$. It preserves the deficit and reverses the sign
$(-1)^{n-|S|}$. All terms cancel. Finally, endpoint rigidity and the
equality of all greedy gaps show that every subset with
$U=V=k>X$ belongs to $\mathcal J$, proving the claim.
\end{proof}

\begin{corollary}[Exceptional top-degree cancellation for $n\equiv2\pmod3$]
\label{cor:both-active-3k-plus-two-top}
Let $n=3r+2\ge11$, and retain the boundary scores
$U(S),V(S),X(S)$ from Lemma~\ref{lem:both-active-boundary-forcing}.
Then
\[
 \sum_{\substack{S\subseteq[n]\\U(S)=V(S)=r+1>X(S)}}
 (-1)^{n-|S|}t^{X(S)}=0.
\]
In particular, if
\[
 B_n^{\mathrm{both}}(t)=
 \sum_{\substack{S\subseteq[n]\\U(S)=V(S)>X(S)}}
 (-1)^{n-|S|}t^{X(S)},
\]
then
\[
 \deg B_n^{\mathrm{both}}(t)\le
 \left\lfloor\frac n3\right\rfloor-1
\]
whenever $B_n^{\mathrm{both}}$ is nonzero and $n\ge11$.
\end{corollary}
\begin{proof}
Write $k=r+1$. Then $n=3k-1$, so the stated cancellation
is exactly Proposition~\ref{prop:both-active-upper-cancellation}.
For the degree bound, if $U(S)=V(S)=\ell>X(S)$, the boundary
choice $U$ supplies an independent $\ell$-set containing $1$
and $n-1$ in the path square. Consecutive vertices in this set
have gaps at least three, hence
$n-2\ge3(\ell-1)$ and
$\ell\le\lfloor(n+1)/3\rfloor$.
Adding either missing wraparound edge can lower the independence
number by at most one; since $U=V>X$, it follows that
$X=\ell-1$. Thus for $n\equiv0,1\pmod3$ the asserted
bound follows immediately, and for $n\equiv2\pmod3$ the only
potentially larger degree is $\ell-1=\lfloor n/3\rfloor$,
whose signed coefficient vanishes by the first assertion.
\end{proof}

\begin{theorem}[Maximal-length both-active cancellation for cycles of length $3k$]
\label{thm:both-active-3k}
Let $n=3k$, $k\ge3$, and use the boundary scores $U(S),V(S),X(S)$
from Lemma~\ref{lem:both-active-boundary-forcing}. Then
\[
 \sum_{\substack{S\subseteq[n]\\ U(S)=V(S)=k>X(S)}}
 (-1)^{n-|S|}t^{\alpha(C_n^2[S])}=-t^{k-1}.
\]
Thus the maximal-length both-active signed contribution is
$-t^{k-1}$.
\end{theorem}
\begin{proof}
For each subset in the sum, Corollary~\ref{cor:both-active-endpoint-rigidity}
gives the left- and right-greedy path-square chains
$L=(a_1<\cdots<a_k)$ and $R=(b_1<\cdots<b_k)$ with
$a_1=1$, $a_k=n-1$, $b_1=2$, and $b_k=n$.
Since the gaps in either chain are at least three and their sum is
$n-2=3(k-1)+1$, each chain has exactly one gap of length four,
all other gaps being three. Write $s$ for the position of the
four-gap of $L$, and $r$ for that of $R$. Endpoint rigidity gives
$2\le s\le k-1$ and $1\le r\le k-2$.

For fixed $(s,r)$, the chains are uniquely determined. Set
$S_0=L\cup R$, $W=A(L)\cap B(R)$, and $F=W\setminus S_0$.
The Boolean-interval lemma says that all subsets with these greedy
chains are precisely $S_0\subseteq S\subseteq W$. Both $L$ and $R$
are independent $k$-sets in $G_n[S]$: their respective endpoint
pairs are $(1,n-1)$ and $(2,n)$, neither of which is an edge of
$G_n$. Thus $U(S)=V(S)=k$ throughout this Boolean interval, and
$\alpha(G_n[S])=k$. Moreover, the only independent $k$-sets of
$C_{3k}^2$ are its three residue classes
\[
 I_0=\{3,6,\ldots,3k\},\quad
 I_1=\{1,4,\ldots,3k-2\},\quad
 I_2=\{2,5,\ldots,3k-1\}.
\]
Consequently, $S$ has deficit one exactly when it contains none
of $I_0,I_1,I_2$; in that case $\alpha(C_n^2[S])=k-1$.

If $r>1$, the first two right-chain entries are $2,5$, and the
first two left-chain entries are $1,4$. Thus $3\in F$.
The right-chain four-gap goes from $b_r=3r-1$ to $b_{r+1}=3r+3$,
so $3r\notin B(R)$ and hence $I_0\nsubseteq W$.
Since $3\notin I_1\cup I_2$, toggling vertex $3$ leaves the
deficit condition unchanged. The alternating sum over this Boolean
interval is therefore zero.

If $s<k-1$, the final gaps of both chains equal three, so
$n-2\in F$. The left-chain four-gap goes from $a_s=3s-2$
to $a_{s+1}=3s+2$, whence $3s+1\notin A(L)$ and
$I_1\nsubseteq W$. As $n-2=3k-2\in I_1$ and
$n-2\notin I_0\cup I_2$, toggling $n-2$ preserves the deficit
condition and reverses the sign. Thus this interval also cancels.

The only pair not covered by these involutions has
$s=k-1$ and $r=1$. Its chains and free set are
\[
 L=(1,4,7,\ldots,3k-5,3k-1),\qquad
 R=(2,6,9,\ldots,3k-3,3k),\qquad
 F=\{5,8,\ldots,3k-4\}.
\]
Here $|S_0|=2k$, $|F|=k-2$, and $I_0,I_1\nsubseteq W$, while
$I_2\subseteq W$ and $I_2\setminus S_0=F$. Hence precisely the
proper subsets $T\subsetneq F$ yield deficit subsets $S=S_0\cup T$.
Their signed contribution is
\[
 (-1)^{3k-2k}\sum_{T\subsetneq F}(-1)^{|T|}t^{k-1}
 =(-1)^k\bigl(-(-1)^{k-2}\bigr)t^{k-1}=-t^{k-1}.
\]
All other greedy-chain intervals cancel, proving the result.
\end{proof}

\begin{lemma}[Coordinatewise comparison of path-square greedy chains]
\label{lem:greedy-offset-bound}
Let $L=(a_1<\cdots<a_k)$ and $R=(b_1<\cdots<b_k)$ be the
left- and right-greedy maximum independent sets of $P_n^2[S]$.
Then $0\le b_j-a_j\le2$ for every $1\le j\le k$.
\end{lemma}
\begin{proof}
The left-greedy algorithm chooses the earliest possible vertex at each
step, so $a_j\le b_j$ by induction. If $b_j\ge a_j+3$, then
$a_1,\ldots,a_j,b_j,\ldots,b_k$ is an independent set of size $k+1$:
its only new consecutive gap is $b_j-a_j\ge3$. This contradicts
$\alpha(P_n^2[S])=k$. Thus $b_j\le a_j+2$.
\end{proof}

\begin{lemma}[Cyclic splicing inequality]\label{lem:cyclic-splicing}
Let $n\ge9$ and suppose $U(S)=V(S)=m>X(S)$.
Write the left- and right-greedy maximum independent sets of
$P_n^2[S]$ as $L=(a_1<\cdots<a_m)$ and
$R=(b_1<\cdots<b_m)$. Then, for every $1\le j<m$,
\[
 a_{j+1}-b_j\le2.
\]
\end{lemma}
\begin{proof}
By Corollary~\ref{cor:both-active-endpoint-rigidity},
$b_1=2$ and $a_m=n-1$. If $a_{j+1}-b_j\ge3$, the
$m$ vertices
\[
 b_1,\ldots,b_j,a_{j+1},\ldots,a_m
\]
are in $S$ and have successive gaps at least three. Their
cyclic wraparound gap is $n+b_1-a_m=3$. Thus they form an
independent $m$-set in $C_n^2[S]$, contradicting
$\alpha(C_n^2[S])=X(S)=m-1$ by
Lemma~\ref{lem:both-active-identifications}.
\end{proof}

\begin{theorem}[Three-or-four-gap restriction for both-active deficits]
\label{thm:five-gap-exclusion}
Let $n\ge9$ and suppose $U(S)=V(S)=m>X(S)$.
For the left- and right-greedy maximum independent sets
$L=(a_1<\cdots<a_m)$ and $R=(b_1<\cdots<b_m)$ of
$P_n^2[S]$, set
\[
 A_j=a_{j+1}-a_j,\quad B_j=b_{j+1}-b_j,\quad
 \delta_j=b_j-a_j.
\]
Then for every $1\le j<m$,
\[
 A_j,B_j\in\{3,4\},\qquad
 A_j=4\Longrightarrow\delta_j=2,\qquad
 B_j=4\Longrightarrow\delta_{j+1}=2.
\]
In particular, when $n=3m+1$ each chain has exactly two
four-gaps, and when $n=3m+2$ each chain has exactly three
four-gaps.
\end{theorem}
\begin{proof}
Put $d_j=a_{j+1}-b_j$. The greedy-chain comparison
(Lemma~\ref{lem:greedy-offset-bound}) gives
$0\le\delta_j\le2$, while the splicing inequality
(Lemma~\ref{lem:cyclic-splicing}) gives $d_j\le2$.
The exact identities
\[
 A_j=d_j+\delta_j,\qquad
 B_j=d_j+\delta_{j+1}
\]
therefore imply $A_j,B_j\le4$. Since both chains are
independent in the path square, $A_j,B_j\ge3$.
If $A_j=4$, the first identity and the bounds
$d_j\le2$, $\delta_j\le2$ force $d_j=\delta_j=2$.
Likewise $B_j=4$ forces $d_j=\delta_{j+1}=2$.

Finally, endpoint rigidity gives
$a_1=1$, $a_m=n-1$, $b_1=2$, and $b_m=n$, so
\[
 \sum_{j=1}^{m-1}(A_j-3)
 =\sum_{j=1}^{m-1}(B_j-3)=n-3m+1.
\]
For $n=3m+1$ this is $2$, and for $n=3m+2$
it is $3$. As each excess is either zero or one,
the asserted counts of four-gaps follow.
\end{proof}

\begin{lemma}[Five local transitions and offset conservation]
\label{lem:five-transitions}
Under the hypotheses of Theorem~\ref{thm:five-gap-exclusion}, put
$d_j=a_{j+1}-b_j$. Then $\delta_j,d_j\in\{1,2\}$, and the only
possible local steps $(\delta_j,\delta_{j+1},d_j;A_j,B_j)$ are
\[
 (1,1,2;3,3),\quad(1,2,2;3,4),\quad(2,1,2;4,3),\quad
 (2,2,1;3,3),\quad(2,2,2;4,4).
\]
Moreover, for $1\le r\le s<m$,
\[
 \sum_{j=r}^s(A_j-B_j)=\delta_r-\delta_{s+1}.
\]
\end{lemma}
\begin{proof}
The gap bounds $A_j,B_j\ge3$ and $d_j,\delta_j,\delta_{j+1}\le2$
force $d_j,\delta_j,\delta_{j+1}\ge1$. Checking the eight possible
triples gives the five displayed rows. Subtraction of
$A_j=d_j+\delta_j$ and $B_j=d_j+\delta_{j+1}$, followed by
summation, proves the telescoping identity.
\end{proof}

\begin{proposition}[Free vertices and excursion counts]
\label{prop:excursion-count}
Assume the hypotheses of Lemma~\ref{lem:five-transitions}.
Let $U_0=L\cup R$, $W=A(L)\cap B(R)$, and $F=W\setminus U_0$.
Then the following union is disjoint:
\[
 F=\{a_j+1:2\le j\le m-1,\ \delta_j=2\}
 \mathbin{\dot\cup}
 \{b_j+1:1\le j<m,\ \delta_j=\delta_{j+1}=1\}.
\]
If $t$ is the number of transitions $1\to2$ in the offset sequence,
and $q$ is the number of steps $(2,2,2;4,4)$, then
\[
 |F|=m-1-t,\quad n=3m-1+t+q,\quad
 |W|=3m-1-t,\quad n-|W|=2t+q.
\]
\end{proposition}
\begin{proof}
Strict interlacing leaves only the candidate points $a_j+1$
when $\delta_j=2$ and $b_j+1$ when $d_j=2$.
The former are covered by both chains, except at the endpoints
where $\delta_1=\delta_m=1$. The latter belong to both coverage
sets precisely for the step $(1,1,2;3,3)$.
Let $z$ be the number of interior indices with offset two. The
numbers of $1\to2$ and $2\to1$ steps both equal $t$, the number
of $2\to2$ steps is $z-t$, and the number of $1\to1$ steps is
$m-1-t-z$. Hence $|F|=z+(m-1-t-z)=m-1-t$.
The left chain has $t+q$ four-gaps, so
$n-2=a_m-a_1=3(m-1)+t+q$. Since $|U_0|=2m$,
the remaining identities follow.
\end{proof}

\begin{proposition}[Minimal-requirement covering identity]
\label{prop:canonical-requirement-cover}
Suppose $U_0\subseteq W\subseteq[n]$ and
$\alpha(C_n^2[U_0])<m$. Put $F=W\setminus U_0$, and let
$\mathcal R$ be the family of inclusion-minimal sets
$J\setminus U_0$ for independent $m$-sets $J\subseteq W$.
Then the signed deficit weight equals
\[
 w=(-1)^{n-|U_0|+|F|}
 \sum_{\substack{\mathcal K\subseteq\mathcal R\\
             \bigcup_{R\in\mathcal K}R=F}}
 (-1)^{|\mathcal K|}.
\]
In particular, if some vertex of $F$ belongs to no minimal
requirement, then $w=0$.
\end{proposition}
\begin{proof}
A set $T\subseteq F$ preserves the deficit if and only if it
contains no member of $\mathcal R$. Inclusion--exclusion over
these forbidden containment events gives
\[
 \sum_{T\subseteq F:\ R\not\subseteq T\ \forall R\in\mathcal R}
 (-1)^{|T|}
 =\sum_{\mathcal K\subseteq\mathcal R}(-1)^{|\mathcal K|}
   \sum_{\bigcup_{R\in\mathcal K}R\subseteq T\subseteq F}(-1)^{|T|}.
\]
The inner sum is zero unless the union equals $F$, in which
case it is $(-1)^{|F|}$. Multiplying by $(-1)^{n-|U_0|}$
proves the formula.
\end{proof}

\begin{theorem}[Uniform canonical-block exhaustiveness without $1\to1$ steps]
\label{thm:general-blocks}
Assume the five-transition hypotheses of
Lemma~\ref{lem:five-transitions}, with
$\delta_1=\delta_m=1$ and no step $(1,1,2;3,3)$.
Let $U_0=L\cup R$, $W=A(L)\cap B(R)$, and $F=W\setminus U_0$.
Partition $F$ into maximal runs $B_1,\ldots,B_k$ of consecutive
vertices at distance three. Then $k=t+q$, every $B_i$ is the
requirement of an independent $m$-set in $W$, and every independent
$m$-set in $W$ contains a complete block. Consequently the
inclusion-minimal requirements are precisely the blocks and
\[
 w=(-1)^t.
\]
\end{theorem}
\begin{proof}
With no $1\to1$ steps, the free-vertex formula of
Proposition~\ref{prop:excursion-count} becomes
$F=\{a_i+1:\delta_i=2\}$, with $2\le i\le m-1$.
The offset-two indices form $t$ nonempty runs. Within each run,
$a_{i+1}-a_i=3$ except at the $q$ special $2\to2$ steps,
where the gap is four. Thus the spacing-three runs number $k=t+q$.

For a block $B=\{a_i+1,\ldots,a_j+1\}$ (consecutive indices), put
\[
 J_B=\{b_1,\ldots,b_{i-1}\}\cup B\cup
       \{a_{j+1},\ldots,a_m\}.
\]
It has $m$ vertices. Its consecutive gaps within each of the three
portions are at least three. At the left junction, either the
preceding step is $1\to2$, with $A_{i-1}=3$ and
$\delta_{i-1}=1$, or it is a special $2\to2$ step,
with $A_{i-1}=4$ and $\delta_{i-1}=2$; in both cases
$(a_i+1)-b_{i-1}=3$. At the right junction, the block
ends at a $2\to1$ step or a special $2\to2$ step,
so $A_j=4$ and $a_{j+1}-(a_j+1)=3$.
Finally $n-a_m+b_1=3$. Thus $J_B$ is an independent
$m$-set and $J_B\setminus U_0=B$.

For exhaustiveness, let $J\subseteq W$ be an independent $m$-set.
The $m$ disjoint intervals $[a_i,a_i+2]$ are cliques in
$C_n^2$ and cover $W$, so $J$ contains exactly one vertex
$c_i=a_i+e_i$ in each. Here $e_i\in\{0,1,2\}$ and,
since $\delta_1=\delta_m=1$ and there are no $1\to1$
steps, $e_1,e_m\in\{0,1\}$. The cyclic closing gap is
$2+e_1-e_m\ge3$, forcing $e_1=1$ and $e_m=0$.
The ordinary gaps give
\[
 e_{i+1}\ge e_i\quad(A_i=3),\qquad
 e_{i+1}\ge e_i-1\quad(A_i=4).
\]
Each canonical block is a consecutive range of indices on which
$e_i$ is nondecreasing. The first block starts at index two and
hence starts with $e_i\ge1$; the last ends at $m-1$, where
$e_i\le1$ because $A_{m-1}=4$ and $e_m=0$.
Between two successive blocks there is either one special
$2\to2$ four-gap, or a $2\to1$ four-gap followed by a
$1\to2$ three-gap. In both cases the first $e$-value of
the next block is at least the last $e$-value of the preceding
block minus one.

If $J$ contains no entire block, then no block has all its
$e$-values equal to one. The first block starts at least at one;
monotonicity forces its last value to be two. Inductively, the
next block starts at least at one, and if it is not entirely
selected its last value is two. This forces the final block to
end at two, contradicting its upper bound one. Thus some entire
block belongs to $J$. Combined with realization, this proves
the assertion about inclusion-minimal requirements.

The blocks partition $F$. Since $\alpha(C_n^2[U_0])=m-1$,
subsets preserving the deficit are exactly those containing no
whole block. Their signed sum factors, giving
\[
 w=(-1)^{n-2m}\prod_{h=1}^k
 \sum_{T\subsetneq B_h}(-1)^{|T|}
 =(-1)^{n-2m+|F|+k}.
\]
Using $n=3m-1+t+q$, $|F|=m-1-t$ and $k=t+q$
shows $w=(-1)^{2m-2+t+2q}=(-1)^t$.
\end{proof}

\begin{lemma}[Ordered forbidden-offset obstruction]\label{lem:ordered-holes}
Let $n=3m-1+K$ with $K\ge1$, let $W\subseteq[n]$, and let $x\in[n]$.
Suppose, after unwrapping the cycle clockwise from $x$, the complement of $W$
contains $K$ vertices
\[
 z_r=x+3a_r+r\qquad(0\le r<K),\qquad
 1\le a_0\le a_1\le\cdots\le a_{K-1}\le m-1.
\]
Then $\alpha(C_n^2[W\setminus N[x]])\le m-2$.
\end{lemma}
\begin{proof}
If $y_1<\cdots<y_{m-1}$ were an independent $(m-1)$-set
in the unwrapped interval $[x+3,x+n-3]$, put
$\sigma_i=y_i-x-3i$. Then
$0\le\sigma_1\le\cdots\le\sigma_{m-1}\le K-1$.
The hole $z_0$ forces $\sigma_{a_0}\ge1$: equality to zero is
forbidden and nonnegativity excludes negative values.
Inductively, if $\sigma_{a_{r-1}}\ge r$, monotonicity gives
$\sigma_{a_r}\ge r$; equality would select the excluded vertex
$z_r$, so $\sigma_{a_r}\ge r+1$. At $r=K-1$ this contradicts
$\sigma_{a_{K-1}}\le K-1$.
\end{proof}

\begin{lemma}[Ordered coverage holes from a $1\to1$ step]
\label{lem:forbidden-offset-holes}
Assume the five-transition hypotheses and put $K=t+q=n-3m+1\ge1$.
If step $j$ is $(1,1,2;3,3)$ and $x=b_j+1$, then the
excluded set $[n]\setminus W$ contains $K$ vertices satisfying
Lemma~\ref{lem:ordered-holes}.
\end{lemma}
\begin{proof}
Write $B_i=b_{i+1}-b_i=3+\varepsilon_i$, where
$\varepsilon_i\in\{0,1\}$, and list the $K$ four-gap indices
$h_1<\cdots<h_K$. At each $h_\ell$, the step is $1\to2$ or
$(2,2,2;4,4)$, so $b_{h_\ell}+1\notin W$.
Indeed, $B_{h_\ell}=4$ forces $d_{h_\ell}=\delta_{h_\ell+1}=2$; if the
step is $1\to2$ then $b_{h_\ell}+1$ lies strictly between $b_{h_\ell}$
and $b_{h_\ell+1}-2=b_{h_\ell}+2$, hence in no right-coverage interval,
while if it is $(2,2,2;4,4)$ then $\delta_{h_\ell}=2$ and
$A_{h_\ell}=4$, so $b_{h_\ell}+1=a_{h_\ell}+3$ lies strictly between
$a_{h_\ell}+2$ and $a_{h_\ell+1}$, hence in no left-coverage interval.
The end of any run of offset twos is a $2\to1$ step and
also gives an excluded vertex $b_u+1$: for such a step
$\delta_u=2$ and $A_u=4$, so $b_u+1=a_u+3$ lies strictly between
$a_u+2$ and $a_{u+1}$, hence outside $A(L)$.
Let $p=\#\{\ell:h_\ell<j\}$. Since $j$ is a $1\to1$
step, $h_p<j<h_{p+1}$ whenever both indices exist.
For $p=0$, choose the holes $b_{h_1}+1,\ldots,b_{h_K}+1$.
Their clockwise distances from $x$ are
\[
 3(h_\ell-j)+(\ell-1),\qquad 1\le\ell\le K.
\]
For $p\ge1$, first choose the $K-p$ holes after $j$, namely
$b_{h_{p+1}}+1,\ldots,b_{h_K}+1$, then after wrapping choose
$b_{h_2}+1,\ldots,b_{h_p}+1$, and finally choose $b_u+1$,
where $u$ is the terminal $2\to1$ step of the offset-two run
containing $h_p$. Such a $u$ satisfies $h_p<u<j$, because
$\delta_j=1$ and every four-gap starts or continues a run of twos.
The distances of the post-$j$ holes are
\[
 3(h_\ell-j)+(\ell-p-1),\qquad p+1\le\ell\le K.
\]
For the wrapped holes they are
\[
 3(m+h_\ell-j)+(K-p+\ell-2),\qquad 2\le\ell\le p,
\]
and the terminal hole has distance
\[
 3(m+u-j)+(K-1).
\]
Here we use that there is no four-gap in $[u,j-1]$: any such index
would be an $h_\ell$ with $h_p<h_\ell<j$, contradicting the choice of
$h_p$ as the largest four-gap below $j$ (and $u$ itself is not a
four-gap index, since $B_u=3$). Hence $b_j-b_u=3(j-u)$, and the
clockwise distance is $n-3(j-u)=3(m+u-j)+(K-1)$.
The number of selected holes is $(K-p)+(p-1)+1=K$ when
$p\ge1$ (and is $K$ when $p=0$).  Their offset terms are
$0,\ldots,K-p-1$ after $j$, then $K-p,\ldots,K-2$
after wrapping, and finally $K-1$ at the terminal hole.
For the terminal hole the index $u$ is strictly larger than $h_p$,
so it comes after every selected wrapped four-gap hole; if $p=1$
there are no wrapped four-gap holes.  Thus the selected holes, in
clockwise order, have exactly the integer offsets $0,1,\ldots,K-1$.
In each case the corresponding multipliers $a_r$ are nondecreasing: the hole indices strictly
increase within each segment, and the transition across the wrap
increases the multiplier from at most $m-1-j$ to at least $m+1-j$.
More explicitly, the post-$j$ multipliers are $h_\ell-j\ge1$,
while the wrapped multipliers are $m+h_\ell-j\le m-1$;
these are ordered since $h_K-j\le m-1-j<m+1-j\le m+h_2-j$.
The terminal multiplier $m+u-j$ is larger than the preceding
wrapped multiplier and at most $m-1$ because $u<j$.
All selected holes therefore lie in $[x+3,x+n-3]$, and their
multipliers satisfy $1\le a_0\le\cdots\le a_{K-1}\le m-1$.
\end{proof}

\begin{theorem}[Uniform signed weights of five-transition intervals]
\label{thm:general-canonical-requirements}
Assume the hypotheses of Lemma~\ref{lem:five-transitions},
with $\delta_1=\delta_m=1$, and let $t$ be the number of
$1\to2$ steps and $q$ the number of special $(2,2,2;4,4)$ steps.
If $K=t+q\ge1$, the signed Boolean-interval weight is
\[
 \boxed{w=\begin{cases}
 (-1)^t,&\text{if no }(1,1,2;3,3)\text{ step occurs},\\
 0,&\text{if at least one such step occurs}.
 \end{cases}}
\]
If $K=0$ and $F\ne\varnothing$, then $w=0$.
\end{theorem}
\begin{proof}
There are two mechanisms.  With no $1\to1$ step, the minimal
requirements are disjoint canonical blocks, so the signed sum
factors into proper-subset sums.  When a $1\to1$ step occurs,
an unused free vertex instead yields a sign-reversing toggle.
If there is no $1\to1$ step, Theorem~\ref{thm:general-blocks}
proves the assertion. Otherwise choose such a step $j$ and
put $x=b_j+1\in F$. If an independent $m$-set $J\subseteq W$
contained $x$, then $J\setminus\{x\}$ would be an independent
$(m-1)$-set in $W\setminus N[x]$. Lemmas~\ref{lem:forbidden-offset-holes}
and~\ref{lem:ordered-holes} rule this out. Thus $x$ belongs
to no minimal requirement; Proposition~\ref{prop:canonical-requirement-cover}
forces $w=0$.
If $K=0$, then $n=3m-1$ and $\alpha(C_n^2)=m-1$;
there are no independent $m$-sets anywhere in the cycle. Therefore
$w=(-1)^{n-2m}(1-1)^{|F|}=0$ when $F$ is nonempty.
\end{proof}

\begin{corollary}[The $3m+2$ signed interval classification]
\label{thm:complete-3m2}
Under the five-transition hypotheses with $n=3m+2$, the signed
interval weight is $(-1)^t$ if there is no $1\to1$ step and
zero otherwise.
\end{corollary}
\begin{proof}
Here $K=n-3m+1=3$, so the assertion is an immediate
specialization of Theorem~\ref{thm:general-canonical-requirements}.
\end{proof}

\begin{theorem}[Complete two-requirement classification for $n=3m+1$]
\label{thm:two-requirement}
Suppose the both-active chains satisfy the five-transition lemma
and $n=3m+1$ (so $t+q=2$). The minimal-requirement family
$\mathcal R$ consists of exactly two nonempty disjoint sets.
More precisely:
\begin{enumerate}
\item If $(t,q)=(1,1)$, let the unique run of offsets two be
$[r,s]$, and let the unique $(2,2,2;4,4)$ step be at $h$,
where $2\le r\le h<s\le m-1$. Then
\[
 \mathcal R=\bigl\{
 \{3j-1:r\le j\le h\},\quad
 \{3j:h+1\le j\le s\}\bigr\}.
\]
The signed interval weight is $-1$ when $r=2,s=m-1$,
and zero otherwise.
\item If $(t,q)=(2,0)$, let the two runs of offsets two be
$[r,s]$ and $[u,v]$, with
$2\le r\le s$, $s+2\le u\le v\le m-1$. Then
\[
 \mathcal R=\bigl\{
 \{3j-1:r\le j\le s\},\quad
 \{3j:u\le j\le v\}\bigr\}.
\]
The signed interval weight is $+1$ when
$r=2,v=m-1,u=s+2$, and zero otherwise.
\end{enumerate}
\end{theorem}
\begin{proof}
For $n=3m+1$, every independent $m$-set has exactly one
cyclic gap of four, and all other gaps are three. Consequently
all such sets are the $n$ distinct sets
\[
 J_a=\{a,a+3,\ldots,a+3(m-1)\}\pmod n.
\]
In particular, for $1\le i\le m$,
\begin{equation}\label{eq:cyclic-start}
 J_{3i}=\{3j-1:1\le j<i\}\cup\{3j:i\le j\le m\}.
\end{equation}
We now derive both the coverage holes and the complete list of
admissible starts, rather than assuming either classification.
A four-gap $A_j=4$ leaves exactly the point $a_j+3$ uncovered
by $A(L)$; a four-gap $B_j=4$ leaves exactly $b_j+1$
uncovered by $B(R)$. The endpoint coverage extends to both
ends of $[n]$, so these are all holes of $W$.

In case (1), the right-chain four-gap at the $1\to2$
transition $j=r-1$ leaves $3r-3$; the simultaneous four-gaps
at $j=h$ leave the same vertex $3h+1$; and the left-chain
four-gap at the $2\to1$ transition $j=s$ leaves $3s+2$.
Thus
\begin{equation}\label{eq:holes-one}
 [n]\setminus W=\{3r-3,3h+1,3s+2\}.
\end{equation}
In case (2), the two $1\to2$ transitions leave the right-chain
holes $3r-3,3u-2$, while the two $2\to1$ transitions leave
the left-chain holes $3s+1,3v+2$. Thus
\begin{equation}\label{eq:holes-two}
 [n]\setminus W=\{3r-3,3s+1,3u-2,3v+2\}.
\end{equation}

To prove that no other $J_a$ lies in $W$, use the bijective
cyclic coordinate $z(x)\equiv-m(x-1)\pmod{3m+1}$.
Adding three to $x$ increases $z(x)$ by one, so every $J_a$
corresponds to a block of $m$ consecutive $z$-positions.
The three residue formulas are
\[
 z(3j)\equiv m+j,\quad z(3j+1)\equiv j,\quad
 z(3j+2)\equiv2m+j+1\pmod{3m+1}.
\]
For (1), the holes in increasing $z$-order are
\[
 h,\quad m+r-1,\quad 2m+s+1.
\]
The three cyclic distances between holes are
\[
 m+r-1-h,\quad m+s-r+2,\quad m-s+h.
\]
The first and third are strictly less than $m$, whereas the
middle one is at least $m+3$. Therefore an allowed block of
$m$ consecutive positions must lie in the middle arc.
Its starting $z$-position ranges from $m+r$ through $m+s+1$.
By $z(3i)=m+i$, the admissible independent sets are exactly
\begin{equation}\label{eq:starts-one}
 J_{3i},\qquad r\le i\le s+1.
\end{equation}
For (2), the holes in increasing $z$-order are
\[
 s,\quad u-1,\quad m+r-1,\quad2m+v+1.
\]
The cyclic distances are
\[
 u-1-s,\quad m+r-u,\quad m+v-r+2,\quad m-v+s.
\]
Only the third distance is at least $m+1$; the other three
are strictly less than $m+1$ and, in fact, at most $m$.
Their allowed interiors have fewer than $m$ positions.
Hence the only possible block starts have $z$-positions
$m+r,\ldots,m+v+1$, giving exactly
\begin{equation}\label{eq:starts-two}
 J_{3i},\qquad r\le i\le v+1.
\end{equation}

In case (1), the chain coordinates and
\eqref{eq:cyclic-start} show that the free vertices required
by $J_{3i}$ are exactly
\[
 R_i=J_{3i}\setminus U_0
 =\{3j-1:r\le j\le h,\ j<i\}
 \cup\{3j:h+1\le j\le s,\ j\ge i\}.
\]
Let $R^+=\{3j-1:r\le j\le h\}$ and
$R^-=\{3j:h+1\le j\le s\}$. Both are nonempty and
disjoint. We have $R_r=R^-$ and $R_{s+1}=R^+$.
If $i\le h+1$, then $R^-\subseteq R_i$; if $i\ge h+1$,
then $R^+\subseteq R_i$. Thus the inclusion-minimal
requirements are exactly $R^-,R^+$.

In case (2), the corresponding requirements are
\[
 R_i=\{3j-1:r\le j\le s,\ j<i\}
 \cup\{3j:u\le j\le v,\ j\ge i\}.
\]
Now $R^+=\{3j-1:r\le j\le s\}$ and
$R^-=\{3j:u\le j\le v\}$, realized at $i=v+1$
and $i=r$, respectively. If $i\le u$ then
$R^-\subseteq R_i$; if $i\ge s+1$ then
$R^+\subseteq R_i$. Since $s+1\le u$, every $i$
satisfies at least one of these conditions. Thus here too
the minimal requirements are precisely $R^-,R^+$.

Finally, since these two requirements are disjoint, the signed
sum over sets containing neither entire requirement vanishes
unless $R^-\cup R^+=F$. The free-vertex formula in
Proposition~\ref{prop:excursion-count} shows that in case
(1) this happens precisely for $r=2,s=m-1$, and in case
(2) precisely for $r=2,v=m-1,u=s+2$: otherwise a free
vertex lies outside both requirements. In the covering case,
the proper-subset sums are $(-1)^{|R^-|+1}$ and
$(-1)^{|R^+|+1}$, respectively, since for every nonempty finite
set $B$ we have
$\sum_{T\subsetneq B}(-1)^{|T|}=(-1)^{|B|+1}$.
Their product is $(-1)^{|R^-|+|R^+|+2}=(-1)^{|F|}$,
so the weight is $(-1)^{n-2m+|F|}=(-1)^t$, giving $-1$ and $+1$,
respectively.
\end{proof}

\begin{theorem}[Pointwise support bounds for both-active deficits]
\label{thm:both-active-support-bounds}
Let $n\ge9$, and suppose $U(S)=V(S)=m>X(S)$, with the
boundary scores of Lemma~\ref{lem:both-active-boundary-forcing}.
Then
\[
 \left\lceil\frac{n+3}{4}\right\rceil
 \le m\le
 \left\lfloor\frac{n+1}{3}\right\rfloor.
\]
In particular, writing $j=X(S)=m-1$, we have
\[
 \left\lfloor\frac{n+2}{4}\right\rfloor
 \le j\le \left\lfloor\frac n3\right\rfloor.
\]
The upper endpoint can occur only when $n\equiv2\pmod3$;
its total signed coefficient vanishes by
Proposition~\ref{prop:both-active-upper-cancellation}.
\end{theorem}
\begin{proof}
By Corollary~\ref{cor:both-active-endpoint-rigidity},
the left-greedy chain has $a_1=1$, $a_2=4$, and $a_m=n-1$.
Theorem~\ref{thm:five-gap-exclusion} says that every successive
gap is either $3$ or $4$.  Its first gap equals $3$, so
\[
 3(m-1)\le n-2=a_m-a_1
 \le3+4(m-2)=4m-5.
\]
Rearranging proves both bounds on $m$.  The degree is
$j=m-1$, and the claimed floor expressions follow.
If $j=\lfloor n/3\rfloor$, the upper bound on $m$
forces $n\equiv2\pmod3$.  The cancellation assertion is
Proposition~\ref{prop:both-active-upper-cancellation},
with its stated range of validity.
\end{proof}

\begin{theorem}[The first staircase coefficient for $n\equiv1\pmod4$]
\label{thm:both-active-first-staircase}
For every integer $m\ge3$, put $n=4m-3$. Then
\[
 \sum_{\substack{S\subseteq[n]\\U(S)=V(S)=m>X(S)}}
 (-1)^{n-|S|}=-1.
\]
Equivalently, the coefficient of $t^{m-1}$ in the
both-active deficit polynomial is $-1$.
\end{theorem}
\begin{proof}
By endpoint rigidity, the first left-greedy gap and the
last right-greedy gap are $3$. The gap restriction gives
$\sum_{j=1}^{m-1}(A_j-3)=(n-2)-3(m-1)=m-2$, and likewise
for the $B_j$; since each gap is $3$ or $4$, each chain
has exactly one gap $3$ and $m-2$ gaps $4$. Hence the two
chains must be
\[
 L=(1,4,8,\ldots,4m-4),\qquad
 R=(2,6,10,\ldots,4m-6,4m-3).
\]
In particular every contributing $S$ contains
\[
 U_0=L\cup R=\{1,n\}\cup\{2,4,6,\ldots,n-1\},
 \qquad |U_0|=2m.
\]
We prove that no further vertex is possible. The vertices
outside $U_0$ are $v=2r+1$, $1\le r\le2m-3$.
For any such $v$, we exhibit an independent $m$-set of
$C_n^2[U_0\cup\{v\}]$.

If $r$ is even, take $v$ together with a maximum set of
nonconsecutive even vertices from
\[
 \{2,4,\ldots,2r-2\}\ \cup\
 \{2r+4,2r+6,\ldots,n-1\}.
\]
The two blocks contain $r-1$ and $2m-r-3$ even vertices,
respectively. Both numbers are odd, and their maximum
nonconsecutive subsets have combined size
\[
 \left\lceil\frac{r-1}{2}\right\rceil+
 \left\lceil\frac{2m-r-3}{2}\right\rceil=m-1.
\]
All chosen vertices are cyclically at distance at least
$3$ from $v$ and from one another, including across the
cyclic boundary.
Here each block has odd length, so its maximum nonconsecutive
subset is unique, namely $2,6,\dots,2r-2$ and
$2r+4,2r+8,\dots,n-1$; the nearest chosen evens to $v$ are
$2r-2$ and $2r+4$, both at distance $3$, and the wrap pair
$2,n-1$ is at cyclic distance $3$.

If $r$ is odd, take $v$ and also take $w=n$ when $r=1$,
or $w=1$ when $r\ge3$. Select a maximum set of
nonconsecutive even vertices from
\[
 \{4,6,\ldots,2r-2\}\ \cup\
 \{2r+4,2r+6,\ldots,n-3\},
\]
interpreting empty blocks as empty. These blocks supply
\[
 \max\!\left(0,\left\lceil\frac{r-2}{2}\right\rceil\right)
 +\max\!\left(0,\left\lceil\frac{2m-r-4}{2}\right\rceil\right)
 =m-2
\]
vertices. Together with $v$ and $w$, they form a cyclically
independent $m$-set. (The endpoint $w$ is at cyclic
distance at least $3$ from $v$ and all the selected evens.)
Each block here has odd length or is empty (empty exactly
for $r=1$ and for $r=2m-3$, where the displayed maxima give
$0$), so the maximum nonconsecutive subset is again unique.
In particular, for $r=2m-3$ the second block is empty, so
$n-3$, which would be at distance $1$ from $v$, is not
selected.
Thus in every case
$\alpha(C_n^2[U_0\cup\{v\}])\ge m$.
A both-active deficit with $U=V=m$ has cyclic
independence number $m-1$, so no such $v$ can belong to $S$.
Therefore $S=U_0$ is the only possible contributing subset.

Conversely, $L$ and $R$ are independent $m$-sets in
$G_n[U_0]$, realizing the two forbidden boundary choices,
so $U(U_0)\ge m$ and $V(U_0)\ge m$. For the reverse
inequalities, a $G_n$-independent set of pattern $(1,n-1)$
contains $1$ and $n-1$, so its remaining vertices lie in
$U_0\cap[4,n-4]$, the $2m-5$ evens $4,\dots,4m-8$, which
form a path of independence number $m-2$; hence
$U(U_0)\le m$. A set of pattern $(2,n)$ contains $2$ and
$n$, so its remaining vertices lie in $U_0\cap[5,n-3]$,
the $2m-5$ evens $6,\dots,4m-6$, again of independence
number $m-2$; hence $V(U_0)\le m$. Thus
$U(U_0)=V(U_0)=m$.
In $C_n^2[U_0]$, the even vertices form a path of
$2m-2$ vertices (in their natural order, consecutive evens
are adjacent; the wrap pair $2,n-1$ is at cyclic distance
$3$), so an independent set using only evens has
size at most $m-1$. An independent set containing $1$
cannot contain $n$, $2$ or $n-1$, the $C_n^2$-neighbours
of $1$ in $U_0$; an independent set containing $n$
cannot contain $1$, $2$ or $n-1$. In either case the
remaining even vertices are $4,\dots,4m-6$, a path of
$2m-4$ vertices of independence number $m-2$, so the set
has size at most $m-1$. Equality is attained by the
alternate evens $\{2,6,10,\dots,4m-6\}$. Hence, by
Lemma~\ref{lem:both-active-identifications},
$\alpha(C_n^2[U_0])=m-1$, so $X(U_0)=m-1$ and
$U(U_0)=V(U_0)=m$.
Finally $n-|U_0|=(4m-3)-2m=2m-3$ is odd, so the
unique contribution is $-1$.
\end{proof}

\subsection{Local neighborhood tests for the next staircase coefficient}
\label{subsec:local-neighborhood-tests}

The following elementary lemmas provide the local neighborhood
calculations used in the adjacent-staircase theorem for $n=4m-4$.
They isolate the extension tests from the subsequent classification
of the forced greedy-chain unions.

\begin{lemma}[Four consecutive neighbors prevent an increase]
\label{lem:four-neighbor-cycle}
Let $H$ be a graph on $2m$ vertices, $m\ge2$, containing a spanning
cycle $v_1v_2\cdots v_{2m}v_1$, and suppose
$\alpha(H)=m-1$. Let $x$ be a new vertex whose neighborhood in $H$
contains four consecutive vertices of this spanning cycle. Then
\[
 \alpha(H+x)=m-1.
\]
\end{lemma}
\begin{proof}
Any independent set containing $x$ avoids those four consecutive
cycle vertices.  The remaining $2m-4$ cycle vertices induce a
spanning path, whose independence number is
$\lceil(2m-4)/2\rceil=m-2$.  Extra edges of $H$ can only
lower this bound.  Hence an independent set containing $x$ has
size at most $1+(m-2)=m-1$.  An independent set avoiding $x$
has size at most $\alpha(H)=m-1$, and equality is already
attained in $H$.
\end{proof}

\begin{lemma}[Three consecutive neighbors: unique candidate]
\label{lem:three-neighbor-cycle}
Let $H$ be a graph on $2m$ vertices containing the spanning cycle
$v_1v_2\cdots v_{2m}v_1$, and let $x$ be a new vertex with
\emph{exactly} the three neighbors
$v_i,v_{i+1},v_{i+2}$ in $H$, with cyclic subscripts. Set
\[
 J_i=\{v_{i+3},v_{i+5},\ldots,v_{i+2m-1}\}.
\]
Then
\[
 \alpha(H-N_H(x))=m-1
 \quad\Longleftrightarrow\quad
 J_i\text{ is independent in }H.
\]
If additionally $\alpha(H)=m-1$, then adjoining $x$ raises the
independence number to $m$ precisely when $J_i$ is independent
in $H$.
\end{lemma}
\begin{proof}
Deleting three consecutive vertices from the spanning cycle
leaves a path on $2m-3$ vertices.  This odd-order path has
independence number $m-1$, and its unique independent set of
that cardinality is its alternating set of vertices in odd
positions, namely $J_i$. Since $H-N_H(x)$ contains that
spanning path, it admits an independent $(m-1)$-set exactly
when $J_i$ survives all additional edges of $H$.
The last assertion follows by adjoining $x$ to $J_i$ and using
$\alpha(H+x)=\max\{\alpha(H),1+\alpha(H-N_H(x))\}$.
\end{proof}

\begin{remark}[Scope of the local tests]
\label{rem:local-tests-scope}
For the $n=4m-4$ family, these lemmas may be applied to a
forced vertex set $T$ once its induced graph
$H=C_n^2[T]$ is shown to contain the indicated spanning cycle
and to satisfy $\alpha(H)=m-1$.  They test individual
extensions, but do not by themselves classify the full
greedy-chain Boolean intervals or prove that their signed
weights alternate.  The following subsection supplies the additional
greedy-chain classification and signed-extension argument needed to
establish $[t^{m-1}]B_{4m-4}^{\mathrm{both}}(t)=-1$; the local
neighbor lemmas alone do not imply that identity.
\end{remark}

\subsection{The adjacent first staircase coefficient}
\label{subsec:adjacent-first-coefficient}

We now complete the local extension analysis for $n=4m-4$.
All indices on the displayed spanning cycles are cyclic unless a
linear range is specified.

\begin{lemma}[Boolean intervals for both-active chains]
\label{lem:both-active-boolean-interval}
Let $S\subseteq[n]$ satisfy $U(S)=V(S)=m>X(S)$, with left- and
right-greedy chains $L=(a_1<\cdots<a_m)$ and
$R=(b_1<\cdots<b_m)$ of $P_n^2[S]$. Put
$A(L)=\bigcup_j[a_j,a_j+2]$, $B(R)=\bigcup_j[b_j-2,b_j]$ and
$W=A(L)\cap B(R)$. Then $L\cup R\subseteq S\subseteq W$, and every
$S'$ with $L\cup R\subseteq S'\subseteq W$ has the same greedy
chains $L,R$; consequently $\alpha(P_n^2[S'])=m$,
$U(S')=V(S')=m$ and $X(S')=\alpha(C_n^2[S'])$.
\end{lemma}
\begin{proof}
Let $x\in S$. Endpoint rigidity gives $a_1=1$ and $a_m=n-1$,
so every vertex of $S$ lies between $a_1$ and $a_m+2$.
If $x\notin A(L)$, it must therefore lie in a gap
$a_j+2<x<a_{j+1}$ for some $j<m$. In particular
$x\ge a_j+3$, so left-greedy on $S$ after picking $a_j$
would pick a vertex $\le x<a_{j+1}$, contradicting its
choice of $a_{j+1}$. Likewise, endpoint rigidity gives
$b_1=2$ and $b_m=n$, so every vertex of $S$ lies between
$b_1-2$ and $b_m$. If $x\notin B(R)$, it must lie in a gap
$b_j<x<b_{j+1}-2$ for some $j<m$. Thus
$x\le b_{j+1}-3$, and right-greedy after choosing
$b_{j+1}$ would choose a vertex $\ge x>b_j$, a
contradiction. Hence $S\subseteq W$. Conversely, let $L\cup R\subseteq S'\subseteq W$.
Left-greedy on $S'$ picks $a_1=1$ first; inductively, after picking
$a_1,\dots,a_j$ the next pick is the smallest $y\in S'$ with
$y\ge a_j+3$. The vertex $a_{j+1}$ qualifies, and any $y<a_{j+1}$
in $W\subseteq A(L)$ lies in some $[a_i,a_i+2]$ with $i\ge j+1$
(as $y\ge a_j+3$), hence $y\ge a_{j+1}$, a contradiction. So the
next pick is $a_{j+1}$. The mirror argument gives the right-greedy
chain $R$. Thus $\alpha(P_n^2[S'])=m$. The chain $L$ is independent
in $G_n[S']$ (it contains $1,n-1$ but not $n$) of size $m$, so
$U(S')=m$, and likewise $R$ gives $V(S')=m$. Finally
$X(S')=\alpha(C_n^2[S'])$: an independent set of $C_n^2[S']$ is
independent in $G_n[S']$ with boundary pattern $\ne(1,n-1),(2,n)$,
and a $G_n[S']$-independent set with any other boundary pattern is
$C_n^2[S']$-independent, since $C_n^2$ adds to $G_n$ only the edges
$\{1,n-1\},\{2,n\}$.
\end{proof}

\begin{lemma}[Consecutive-pair deletion]
\label{lem:consecutive-pair-deletion}
Let $H_m$ be the cycle $v_1\cdots v_{2m}v_1$ with the
additional edges $v_1v_{2m-1}$ and $v_2v_{2m}$, where $m\ge3$.
Then, for every cyclically consecutive pair,
\[
 \alpha(H_m-\{v_i,v_{i+1}\})=m-1.
\]
\end{lemma}
\begin{proof}
The remaining spanning path has $2m-2$ vertices, so its
independence number is at most $m-1$. For odd
$1\le i\le2m-1$, choose
\[
 \{v_2,v_4,\ldots,v_{i-1}\}\cup
 \{v_{i+2},v_{i+4},\ldots,v_{2m-1}\}.
\]
For even $2\le i\le2m-2$, choose
\[
 \{v_3,v_5,\ldots,v_{i-1}\}\cup
 \{v_{i+2},v_{i+4},\ldots,v_{2m}\}.
\]
Both sets have $m-1$ vertices and avoid both extra edges:
the former omits $v_1,v_{2m}$ and the latter omits $v_1,v_2$.
For the cyclic pair $\{v_{2m},v_1\}$ choose
$\{v_2,v_4,\ldots,v_{2m-2}\}$. These are independent
sets of the required cardinality.
\end{proof}

\begin{lemma}[Greedy intervals and their free vertices]
\label{lem:adjacent-free-vertices}
Let $m\ge4$, $n=4m-4$, and $0\le r\le2m-6$. Put
\[
 T_r=\{1,n\}\cup\{2,4,\ldots,2r+4\}
       \cup\{2r+5,2r+7,\ldots,n-1\}.
\]
Every both-active deficit with $U=V=m$ belongs to one of the
$2m-5$ greedy-chain Boolean intervals with lower endpoint $T_r$.
Writing $T_r=\{v_1<\cdots<v_{2m}\}$, we have
\[
 v_i=\begin{cases}
 1&i=1,\\
 2i-2&2\le i\le r+3,\\
 2i-3&r+4\le i\le2m-1,\\
 n&i=2m.
 \end{cases}
\]
The free vertices $F_r$ of these intervals are exactly
\[
\begin{aligned}
 F_0&=\{3\}\cup\{4j:2\le j\le m-2\},\\
 F_{2q}&=\{4j+1:1\le j\le q\}
       \cup\{4j:q+2\le j\le m-2\},
       &&1\le q\le m-4,\\
 F_{2q+1}&=\{4j+1:1\le j\le q+1\}
       \cup\{4j:q+2\le j\le m-2\},
       &&0\le q\le m-4,\\
 F_{2m-6}&=\{4j+1:1\le j\le m-3\}\cup\{n-2\}.
\end{aligned}
\]
Moreover, the neighbors in $T_r$ of any ordinary free vertex
are two consecutive $v_i$; the exceptions are $3$ in $F_0$
and $n-2$ in $F_{2m-6}$ (four consecutive neighbors),
and, for odd $r=2q+1$, the vertices
$x=4q+5$, $y=4q+8$, whose respective neighborhoods are
\[
 \{v_{2q+3},v_{2q+4},v_{2q+5}\},\qquad
 \{v_{2q+4},v_{2q+5},v_{2q+6}\}.
\]
\end{lemma}
\begin{proof}
The endpoint rigidity and gap restrictions established earlier
force exactly two gaps of length three in each greedy chain.
The first left gap and last right gap are already three.
Let their other positions be $s$ for the left chain and $t$
for the right chain. We justify the relative-position restriction
used below. Set $\delta_j=b_j-a_j$ and
$d_j=a_{j+1}-b_j$. The previous lemmas give
$0\le\delta_j\le2$ and $d_j\le2$, hence
\[
 A_j=4\Longrightarrow\delta_j=2,\qquad
 B_j=4\Longrightarrow\delta_{j+1}=2.
\]
Moreover $B_j-A_j=\delta_{j+1}-\delta_j$.
If $s=t$, the required restriction is immediate. Suppose
$s\ne t$. First consider $s\le m-2$. Then $A_s=3$ and
$B_s=4$, because the only right three-gap positions are
$t$ and $m-1$. Consequently $\delta_{s+1}=2$ and
$\delta_s=1$. Since $s\ge2$, the latter equality forces
$B_{s-1}=3$: otherwise $B_{s-1}=4$ would imply
$\delta_s=2$. As $s-1\le m-3$, the only possible right
three-gap position there is $t$, so $t=s-1$.

It remains to handle the endpoint $s=m-1$, where the preceding
argument cannot be used because $B_{m-1}=3$. Endpoint rigidity
gives $\delta_m=1$. Both final gaps have length three, hence
$B_{m-1}-A_{m-1}=0$ and therefore
$\delta_{m-1}=\delta_m=1$. If $t\ne m-2$, then
$B_{m-2}=4$, which would force $\delta_{m-1}=2$, a
contradiction. Thus $t=m-2=s-1$. In all cases,
\[
 t\le s\le t+1.
\]
The possibilities $(s,t)=(q+2,q+1)$
and $(q+2,q+2)$ give respectively $r=2q$ and $r=2q+1$.
Indeed, $A_j=3$ iff $j\in\{1,s\}$ and $B_j=3$ iff
$j\in\{t,m-1\}$, so
\[
a_1=1,\quad a_j=4j-4-\mathbf 1_{j>s}\ (j\ge2),\qquad
b_1=2,\quad b_j=4j-2-\mathbf 1_{j>t}-\mathbf 1_{j>m-1}\ (j\ge2),
\]
with $a_m=n-1$ and $b_m=n$. If $(s,t)=(q+2,q+1)$ then
$a_j=4j-4$ for $2\le j\le q+2$ and $a_j=4j-5$ for $j\ge q+3$,
while $b_j=4j-2$ for $2\le j\le q+1$ and $b_j=4j-3$ for
$q+2\le j\le m-1$; the union is the evens $2,4,\dots,4q+4$,
the odds $1$ and $4q+5,4q+7,\dots,4m-5$, and $n$, i.e.\ $T_{2q}$.
If $(s,t)=(q+2,q+2)$ then $L$ is as above and
$b_j=4j-2$ for $2\le j\le q+2$, $b_j=4j-3$ for $j\ge q+3$; the
union is the evens $2,4,\dots,4q+6$, the odds $1$ and
$4q+7,4q+9,\dots,4m-5$, and $n$, i.e.\ $T_{2q+1}$. The sorted
list is $v_1=1$, $v_i=2i-2$ for $2\le i\le r+3$,
$v_i=2i-3$ for $r+4\le i\le2m-1$, $v_{2m}=n$.
In the even case $r=2q$, the offsets satisfy $\delta_1=\delta_m=1$,
$\delta_{q+2}=1$, and all other interior offsets equal $2$;
all splicing gaps $d_j$ equal $2$. In the odd case $r=2q+1$,
$\delta_1=\delta_m=1$ and all interior offsets equal $2$,
while $d_j=2$ except for $d_{q+2}=1$. These follow by
subtracting the displayed formulas for $a_j$ and $b_j$.

For completeness, the free-set calculation uses the exact
left- and right-greedy coverage conditions
\[
 x\in W_r\ \Longleftrightarrow\quad
 \bigl(\exists a\in L_r:\ a\le x\le a+2\bigr)
 \text{ and }
 \bigl(\exists b\in R_r:\ b-2\le x\le b\bigr),
\]
i.e.\ $W_r=A(L_r)\cap B(R_r)$. The windows $[a_j,a_j+2]$ abut
iff $A_j=3$ and otherwise leave the single gap $\{a_j+3\}$, so
$[n]\setminus A(L_r)=\{a_j+3:A_j=4,\ j<m\}$; likewise
$[n]\setminus B(R_r)=\{b_j+1:B_j=4,\ j<m\}$. Let
$x\in[n]\setminus T_r$ lie in $[a_j,a_j+2]$. If $x=a_j+1$ then
$x\in B(R_r)$ always, since
$[b_j-2,b_j]=[a_j+\delta_j-2,a_j+\delta_j]\ni a_j+1$ for
$\delta_j\in\{1,2\}$, and $x\in T_r$ iff $\delta_j=1$ (then
$x=b_j$); so $x$ is free iff $\delta_j=2$. From the offset table
in the proof above this gives $x=4j+1$ with $1\le j\le q$ for
$r=2q$ and $1\le j\le q+1$ for $r=2q+1$ on the left, and
$x=4j$ with $q+2\le j\le m-2$ on the right, in either case.
If $x=a_j+2$ with $\delta_j=1$, then $x=b_j+1$ and
$x\in B(R_r)$ iff $B_j=3$. Here $\delta_j=1$ at $j=1$ and
$j=q+2$ (case $r=2q$) and at $j=1$ (case $r=2q+1$). In case
$r=2q$, $B_1=2+\delta_2$, so $B_1=3$ iff $q=0$, giving the free
vertex $3$ at $r=0$ only; and $B_{q+2}=2+\delta_{q+3}=4$ unless
$q=m-3$, where $B_{m-1}=3$ gives the free vertex $a_{m-1}+2=n-2$
at $r=2m-6$ only. In case $r=2q+1$, $B_1=4$, so $3\notin W_r$.
This proves all four free-set formulas.

Finally, a vertex $x$ of the ambient cycle square is adjacent
to exactly those vertices of $T_r$ lying in
$\{x-2,x-1,x+1,x+2\}$, interpreted modulo $n$.
Let $x=a_j+1$ with $\delta_j=2$. Its possible $T_r$-neighbours
are $\{a_j-1,a_j,a_j+2,a_j+3\}$; here $a_j\in L_r$,
$a_j+2=b_j\in R_r$, $a_j-1\in T_r$ iff $d_{j-1}=1$, and
$a_j+3\in T_r$ iff $A_j=3$. Hence
$|N(x)\cap T_r|=2+\mathbf 1_{d_{j-1}=1}+\mathbf 1_{A_j=3}$.
In case $r=2q$ we have $d\equiv2$ and $A_j=3$ only at $j=s$,
where $\delta_s=1$, so $a_s+1$ is not free; every ordinary free
vertex has exactly the two consecutive neighbours
$\{a_j,b_j\}$. In case $r=2q+1$, $A_j=3$ at $j=s$ gives
$x=a_s+1=4q+5$ with neighbours $\{a_s,b_s,a_{s+1}\}$, and
$d_{j-1}=1$ only at $j=s+1$ gives $y=a_{s+1}+1=4q+8$ with
neighbours $\{b_s,a_{s+1},b_{s+1}\}$; counting elements
$\le a_s=4q+4$ in $T_r$ gives $a_s=v_{2q+3}$, $b_s=v_{2q+4}$,
$a_{s+1}=v_{2q+5}$, $b_{s+1}=v_{2q+6}$, i.e.\ the two displayed
triples of consecutive $v$'s. Every other free vertex has
exactly two consecutive neighbours. At the extreme endpoints,
$N(3)=\{1,2,4,5\}=\{v_1,v_2,v_3,v_4\}$ and
$N(n-2)=\{n-4,n-3,n-1,n\}=\{v_{2m-3},v_{2m-2},v_{2m-1},v_{2m}\}$,
four consecutive neighbours in each case.
\end{proof}

\begin{theorem}[The adjacent first staircase coefficient]
\label{thm:both-active-adjacent-first}
For every integer $m\ge4$, with $n=4m-4$,
\[
 \sum_{\substack{S\subseteq[n]\\U(S)=V(S)=m>X(S)}}
 (-1)^{n-|S|}=-1.
\]
Equivalently, $[t^{m-1}]B_n^{\mathrm{both}}(t)=-1$.
\end{theorem}
\begin{proof}
Fix $r$. The $v$-gaps are $g_1=1$, $g_i=2$ for $2\le i\le r+2$,
$g_{r+3}=1$, $g_i=2$ for $r+4\le i\le2m-2$, $g_{2m-1}=1$ and
$g_{2m}=1$ (wraparound $n\to1$). Two consecutive gaps sum to $2$
only for the pairs $(g_{2m-1},g_{2m})$ and $(g_{2m},g_1)$, and
any three consecutive gaps sum to at least $3$; hence two
distinct $v_i,v_j$ are at cyclic distance $\le2$ in $[n]$ iff
they are consecutive in the $v$-cycle or
$\{i,j\}=\{1,2m-1\}$ or $\{2,2m\}$. So the induced graph
$C_n^2[T_r]$ is the $2m$-cycle in the ordered vertices $v_i$
with two extra edges $v_1v_{2m-1}$ and $v_2v_{2m}$, i.e.\ the
graph $H_m$ of Lemma~\ref{lem:consecutive-pair-deletion}. Its
independence number is $m-1$: the only two independent $m$-sets
of the even cycle are its alternating classes, and each violates
one of these extra edges. An alternating class with one
vertex removed supplies an independent $(m-1)$-set.
The greedy chains of $P_n^2[T_r]$ are $L$ and $R$: left-greedy
picks $a_{j+1}$ as the smallest element of $T_r$ at distance
$\ge3$ from $a_j$, since the only element of $T_r$ strictly
between $a_j$ and $a_{j+1}$ is $b_j$, at distance $\delta_j\le2$,
and similarly on the right. The chain $L$ contains $1,n-1$ and
$R$ contains $2,n$; both are independent in $G_n[T_r]$, and
$\alpha(P_n^2[T_r])=m$, so $U(T_r)=V(T_r)=m$. With
$\alpha(C_n^2[T_r])=m-1$, $T_r$ itself is a both-active deficit.

By Lemma~\ref{lem:both-active-boolean-interval}, every set $S$ in the interval $[T_r,W_r]$ has $U(S)=V(S)=m$, and hence contributes to the both-active deficit sum precisely when $\alpha(C_n^2[S])=m-1$.

By Lemma~\ref{lem:adjacent-free-vertices}, all free vertices
except the listed exceptions have two consecutive neighbors
in $T_r$. Lemma~\ref{lem:consecutive-pair-deletion} shows that
adjoining any one of them creates a cyclic independent
$m$-set, so none can occur in a deficit extension.
The four-neighbor lemma (Lemma~\ref{lem:four-neighbor-cycle})
allows the single exceptional extension at each endpoint
$r=0,2m-6$.

For odd $r=2q+1$, the two exceptional vertices are
$x=4q+5$ and $y=4q+8$. By
Lemma~\ref{lem:three-neighbor-cycle}, adjoining either
one preserves independence number $m-1$: the unique
alternating $(m-1)$-set remaining after deleting its three
neighbors contains respectively both ends of the extra
edge $v_2v_{2m}$ or $v_1v_{2m-1}$. Indeed, for $x$ the
neighbour index is $i=2q+3$ and
$J_i=\{v_{2q+6},v_{2q+8},\dots,v_{2m},v_2,v_4,\dots,v_{2q+2}\}$
(indices mod $2m$) contains $v_2$ and $v_{2m}$; for $y$ the
neighbour index is $i=2q+4$ and
$J_i=\{v_{2q+7},v_{2q+9},\dots,v_{2m-1},v_1,v_3,\dots,v_{2q+3}\}$
contains $v_1$ and $v_{2m-1}$. Hence in each case $J_i$ is
not independent in $H$, $\alpha(H-N(x))=m-2$ (resp.\ for $y$),
and adjoining the vertex keeps the independence number $m-1$.
They cannot be adjoined together. Indeed,
\[
 I=\{v_2,v_4,\ldots,v_{2q+2}\}
 \cup\{v_{2q+7},v_{2q+9},\ldots,v_{2m-1}\}
\]
has $(q+1)+(m-q-3)=m-2$ vertices and is independent in $H$:
its elements are $2$ cycle-steps apart except for the jump
from $v_{2q+2}$ to $v_{2q+7}$ and the wraparound from
$v_{2m-1}$ to $v_2$, which are $5$ resp.\ $3$ cycle-steps; it
contains $v_2$ but not $v_{2m}$, and $v_{2m-1}$ but not $v_1$,
so neither extra edge is violated. It avoids both
neighborhoods, since
$N(x)=\{v_{2q+3},v_{2q+4},v_{2q+5}\}$ and
$N(y)=\{v_{2q+4},v_{2q+5},v_{2q+6}\}$ are disjoint from $I$.
Since $y-x=3$, $x$ and $y$ are non-adjacent in $C_n^2$ and
non-adjacent to $I$ (any $C_n^2$-neighbour of $x$ in $T_r$ lies
in $N(x)$, and likewise for $y$), so $I\cup\{x,y\}$ is a cyclic
independent $m$-set. Thus the permissible extensions
are $\varnothing$ at an interior even $r$,
$\varnothing,\{x\},\{y\}$ at an odd $r$, and
$\varnothing$ together with the one exceptional singleton
at either endpoint.

Every $T_r$ has $2m$ vertices, so its sign is
$(-1)^{n-2m}=1$, while a one-vertex extension has sign
$-1$. The endpoint intervals have weight $0$, the
$m-4$ interior even intervals have weight $+1$, and the
$m-3$ odd intervals have weight $-1$. Summing gives
$(m-4)-(m-3)=-1$.
\end{proof}

\begin{theorem}[Maximal-length both-active cancellation for cycles of length $3k+1$]
\label{thm:both-active-3k-plus-one}
Let $k\ge4$, $n=3k+1$, and let $U(S),V(S),X(S)$ be the
boundary-choice scores defined above. Then
\[
 \sum_{\substack{S\subseteq[n]\\U(S)=V(S)=k>X(S)}}
 (-1)^{n-|S|}t^{k-1}=-t^{k-1}.
\]
This is an identity for the maximal-length both-active part of the
deficit, not for the entire deficit polynomial.
\end{theorem}
\begin{proof}
By Corollary~\ref{cor:both-active-endpoint-rigidity} and
Theorem~\ref{thm:five-gap-exclusion}, the left chain $L=(a_j)$
has exactly two gaps of four, at positions $s<t$, and the right
chain $R=(b_j)$ has exactly two, at positions $p<q$. Thus
\[
 2\le s<t\le k-1,\qquad 1\le p<q\le k-2,
\]
and
\begin{equation}\label{eq:four-gap-chains}
 a_j=3j-2+\mathbf{1}_{j>s}+\mathbf{1}_{j>t},\qquad
 b_j=3j-1+\mathbf{1}_{j>p}+\mathbf{1}_{j>q}.
\end{equation}
We work separately on the Boolean interval belonging to each
pair of greedy chains. Write
\[
 U_0=L\cup R,\qquad W=A(L)\cap B(R),\qquad F=W\setminus U_0.
\]
The sets with these greedy chains are exactly the interval
$U_0\subseteq S\subseteq W$: if $x\in S$ and $x\notin A(L)$,
greedy minimality is contradicted by $x$ (or $x$ together with
$L$ gives an independent $(k+1)$-set), and symmetrically for
$B(R)$; conversely the windows $[a_j,a_j+2]$ are disjoint and
contain $W$, so for $U_0\subseteq S'\subseteq W$ the greedy
chains of $S'$ are again $L,R$. Every $S$ in this interval has
$\alpha(P_n^2[S])=k$, since $L\subseteq S$ and each window
$[a_j,a_j+2]$ contains at most one vertex of any $P_n^2$-independent
set.
Both $L$ and $R$ are independent $k$-sets in $G_n[S]$:
$L$ uses $1,n-1$, while $R$ uses $2,n$.
Consequently $U(S)=V(S)=k$ throughout the interval. By
Lemma~\ref{lem:both-active-identifications} and the cyclically independent
$(k-1)$-set $L\setminus\{1\}$, we have
$\alpha(C_n^2[S])\in\{k-1,k\}$ on the interval, and its members
contribute precisely when $\alpha(C_n^2[S])=k-1$.

The maximum independent sets of $C_{3k+1}^2$ are the $3k+1$
sets $J_r=\{r,r+3,\ldots,r+3(k-1)\}$, with cyclic indices.
A four-gap of $L$ leaves one vertex outside $A(L)$; similarly
for $R$. Formula~\eqref{eq:four-gap-chains} gives
\[
 [n]\setminus W=\{3s+1,\,3t+2,\,3p,\,3q+1\},
\]
as a set, with $3s+1=3q+1$ when $s=q$ (the other values are
distinct modulo $3$), so $|[n]\setminus W|$ is $3$ if $s=q$ and
$4$ otherwise.
Order the vertices cyclically by repeated addition of $3$.
The four missing positions in this order, using coordinate
$z(v)\equiv-k(v-1)\pmod{3k+1}$, are
\[
 s,\quad q,\quad k+p,\quad 2k+t+1.
\]
Put $u=\min(s,q)$ and $v=\max(s,q)$.  The cyclically
ordered excluded positions are $u,v,k+p,2k+t+1$; when
$s=q$, the first two positions coincide and are counted only
once.  The successive distances, including the zero distance
when $u=v$, are
\[
 v-u,\qquad k+p-v,\qquad k+t-p+1,\qquad k-t+u.
\]
The first is less than $k$ because $s,q\in[1,k-1]$.
The second is less than $k$ since $v\ge q>p$, and the
fourth is at most $k-1$ since $u\le s<t$.  A cyclic arc
between excluded positions at distance $d$ has exactly $d-1$
allowed positions, so an independent $k$-set, which is a block
of $k$ consecutive positions in the $z$-coordinate, cannot
lie in any of these three arcs.  The remaining arc runs from
$k+p$ to $2k+t+1$ and admits precisely the block-start
coordinates
\[
 k+p+1\le a\le k+t+1.
\]
Indeed, its $d-1=k+t-p$ allowed positions contain a block
of $k$ positions exactly when its start lies in that range.
Since $z(3j)\equiv k+j\pmod{3k+1}$, the maximum
independent sets contained in $W$ are exactly
\begin{equation}\label{eq:surviving-J}
 J_{3j},\qquad p+1\le j\le t+1,
\end{equation}
where an empty index range means there are none. For $1\le j\le k$,
\[
 J_{3j}=\{3h-1:1\le h<j\}\cup\{3h:j\le h\le k\}.
\]
The piecewise forms of the two chains are
\[
 a_h=\begin{cases}
 3h-2,&h\le s,\\3h-1,&s<h\le t,\\3h,&h>t,
 \end{cases}
 \qquad
 b_h=\begin{cases}
 3h-1,&h\le p,\\3h,&p<h\le q,\\3h+1,&h>q.
 \end{cases}
\]
Thus the additional free vertices required to include $J_{3j}$
are exactly
\begin{equation}\label{eq:requirements}
 M_j=J_{3j}\setminus U_0
 =\{3h-1:p<h<j,\ h\le s\}
 \cup\{3h:\max(j,q+1)\le h\le t\}.
\end{equation}
Their union is
\[
 P\cup Q,\qquad
 P=\{3h-1:p<h\le s\},\qquad
 Q=\{3h:q<h\le t\}.
\]
If a free vertex is outside $P\cup Q$, toggling that vertex
preserves every containment condition $J_{3j}\subseteq S$ and
reverses the sign. Hence the interval has zero signed contribution.

We now determine the pairs that can escape this cancellation.
If $p>1$, vertex $3$ is free and is outside $P\cup Q$;
therefore $p=1$ is necessary. If $t<k-1$, vertex $n-2$
is free and outside $P\cup Q$; therefore $t=k-1$ is
necessary. Lemma~\ref{lem:greedy-offset-bound} forces $q\ge s$:
if $q<s$, the offset at index $q+1$ would equal $3$.
If $q\ge s+2$, vertex $3s+4$ is free but lies outside
$P\cup Q$, so it also gives cancellation. Thus a nonzero
interval must have
\[
 p=1,\qquad t=k-1,\qquad q=s\ \text{or}\ q=s+1.
\]

Write $r=s$. These are the two families
\[
 \mathrm A_r:(s,t;p,q)=(r,k-1;1,r),\quad 2\le r\le k-2,
\]
\[
 \mathrm B_r:(s,t;p,q)=(r,k-1;1,r+1),\quad 2\le r\le k-3.
\]
In either family, the offsets $b_j-a_j$ are $1$ or $2$ and the
splicing gaps $d_j=a_{j+1}-b_j$ are $1$ or $2$, so
$a_j<b_j<a_{j+1}$ for $j<k$; hence $L\cap R=\varnothing$ and
$|U_0|=2k$.
For family $\mathrm A_r$, the two distinct coverage holes at
$3r+1$ coincide; hence $|W|=n-3$, $|F|=k-2$, and
\[
 F=P_r\sqcup Q_r,\quad
 P_r=\{5,8,\ldots,3r-1\},\quad
 Q_r=\{3r+3,3r+6,\ldots,3k-3\}.
\]
For family $\mathrm B_r$ the four holes are distinct, so
$|W|=n-4$, $|F|=k-3$, and
\[
 F=P_r\sqcup Q'_r,\qquad
 Q'_r=\{3r+6,3r+9,\ldots,3k-3\}.
\]
All displayed blocks are nonempty in their indicated ranges.
We verify the minimal requirements explicitly.  In both
families, the surviving indices are $2\le j\le k$.
For $\mathrm A_r$ the requirement formula becomes
\[
 M_j=\{3h-1:2\le h\le r,\ h<j\}
 \cup\{3h:r<h\le k-1,\ h\ge j\}.
\]
For $\mathrm B_r$ the same formula holds with $r<h$
replaced by $r+1<h$ in the second set.  Thus for
$2\le j\le r$ the second block is contained in $M_j$,
whereas for $r+1\le j\le k$ the first block is
contained in $M_j$.  Moreover, $M_2$ is exactly the second
block and $M_k=P_r$.  Hence a subset $T\subseteq F$
contains some requirement $M_j$ if and only if it contains
$P_r$ or the entire second block.  The two blocks are
nonempty and disjoint, so $S=U_0\cup T$ is a deficit set
exactly when $T$ contains neither block in full.
For disjoint nonempty blocks $P,Q$ the signed sum of such $T$ is
\[
 \left(\sum_{T_P\subsetneq P}(-1)^{|T_P|}\right)
 \left(\sum_{T_Q\subsetneq Q}(-1)^{|T_Q|}\right)
 =(-1)^{|P|+|Q|}.
\]
Multiplying by $(-1)^{n-|U_0|}=(-1)^{k+1}$, each
$\mathrm A_r$ contributes $-t^{k-1}$ and each
$\mathrm B_r$ contributes $+t^{k-1}$.
There are $k-3$ configurations of type $\mathrm A$ and $k-4$
of type $\mathrm B$, so the total is
\[
 -(k-3)t^{k-1}+(k-4)t^{k-1}=-t^{k-1},
\]
as asserted.
\end{proof}

\subsection{Boundary cancellation for two additional short gaps}
\label{subsec:two-gap-boundary}

The following result isolates a uniform boundary cancellation in the
case $n=4m-5$.  It is proved by examining the minimal
cyclic requirements in the corresponding boundary intervals.

Fix $m\ge5$, put $n=4m-5$, and choose two-element gap-position sets
$P=\{p_1<p_2\}\subseteq\{2,\ldots,m-1\}$ and
$Q=\{q_1<q_2\}\subseteq\{1,\ldots,m-2\}$ satisfying
$q_j\le p_j\le q_j+1$.  Let $L_H,R_H$ be the associated
left- and right-greedy chains with their prescribed endpoint
three-gaps, where $H=\{h_1,h_2\}$ and $h_j=p_j+q_j$.
Set
\[
 U_H=L_H\cup R_H,\qquad
 W_H=\left(\bigcup_{a\in L_H}[a,a+2]\right)
 \cap\left(\bigcup_{b\in R_H}[b-2,b]\right)\cap[n],
 \qquad F_H=W_H\setminus U_H.
\]
For $J\subseteq W_H$ an independent $m$-set of $C_n^2$, its
\emph{requirement} is $M_J=J\setminus U_H$.  The signed interval
weight is
\[
 w_m(H)=\sum_{\substack{T\subseteq F_H\\
 \alpha(C_n^2[U_H\cup T])=m-1}}
 (-1)^{n-|U_H|-|T|}.
\]

\begin{lemma}[Explicit offsets and free vertices for two additional short gaps]
\label{lem:two-gap-free-vertices}
Fix the two-gap data of Subsection~\ref{subsec:two-gap-boundary}, and
write $L_H=(a_1,\ldots,a_m)$ and $R_H=(b_1,\ldots,b_m)$.
Put $\delta_i=b_i-a_i$ and $A_i=a_{i+1}-a_i$ for $i<m$.
Then
\begin{align*}
 a_i&=4i-3-\mathbf1_{\{i>1\}}-
       \sum_{j=1}^2\mathbf1_{\{i>p_j\}},\\
 b_i&=4i-2-\mathbf1_{\{i>m-1\}}-
       \sum_{j=1}^2\mathbf1_{\{i>q_j\}}.
\end{align*}
In particular, $\delta_1=\delta_m=1$, and for $2\le i\le m-1$,
\[
 \delta_i=2-\sum_{j=1}^2\mathbf1_{\{q_j<i\le p_j\}}.
\]
Consequently, $\delta_i\in\{1,2\}$, and the free-vertex set has the
disjoint description
\begin{equation}\label{eq:two-gap-free-set}
 F_H=\{a_i+1:\delta_i=2\}\ \sqcup\
 \{b_i+1:1\le i<m,\ \delta_i=\delta_{i+1}=1,\ A_i=3\}.
\end{equation}
If $e_j=p_j-q_j\in\{0,1\}$, then
\begin{equation}\label{eq:two-gap-free-count}
\boxed{\begin{aligned}
 |F_H|={}&m-2-e_1-e_2
 +\mathbf1_{\{q_1=1\}}\\
 &+\mathbf1_{\{e_1=1,\ q_2=p_1\}}
 +\mathbf1_{\{p_2=m-1,\ q_2=m-2\}}.
\end{aligned}}
\end{equation}
\end{lemma}
\begin{proof}
The left chain starts at $1$ and has three-step gaps precisely
at indices $1,p_1,p_2$; all its other gaps have length four.
The right chain starts at $2$ and has three-step gaps precisely
at $q_1,q_2,m-1$. Summing these successive gaps proves the
coordinate formulas. Subtracting them proves the offset formula.
Since $p_j-q_j\in\{0,1\}$ and the two index pairs are ordered,
the index intervals $(q_j,p_j]$ are disjoint. Thus the only
interior indices with $\delta_i=1$ are $i=p_j$ for those $j$
with $e_j=1$.

The chain vertices interlace as
$a_i<b_i<a_{i+1}$ for $i<m$. Every vertex of
$W_H\setminus U_H$ therefore lies either strictly between
$a_i$ and $b_i$, or strictly between $b_i$ and $a_{i+1}$.
In the first interval, $\delta_i\le2$, so the sole possible
vertex is $a_i+1$, occurring exactly when $\delta_i=2$;
it belongs to both coverage windows.
In the second interval, the interlacing bound
$a_{i+1}-b_i\le2$ leaves only the candidate $b_i+1$.
For it to occur we must have $a_{i+1}-b_i=2$;
its membership in the left coverage region requires
$b_i+1\le a_i+2$, hence $\delta_i=1$;
its membership in the right coverage region requires
$b_i+1\ge b_{i+1}-2$, hence $\delta_{i+1}=1$.
These conditions are equivalent to
$\delta_i=\delta_{i+1}=1$ and $A_i=3$.
This proves the disjoint union in~\eqref{eq:two-gap-free-set}.

There are $m-2-e_1-e_2$ interior indices with offset two.
For the second type, $A_i=3$ occurs only at $i=1,p_1,p_2$.
At $i=1$ both offsets are one exactly when $q_1=1$.
At $i=p_1$ they are one exactly when
$e_1=1$ and $q_2=p_1$.
At $i=p_2$ they are one exactly when
$p_2=m-1$ and $q_2=m-2$.
Adding these contributions gives~\eqref{eq:two-gap-free-count}.
\end{proof}

\begin{lemma}[Exchange at the left boundary]
\label{lem:two-gap-boundary-exchange}
If $h_1=3$, then $3\in F_H$ and $2,4\in U_H$.
For any independent $m$-set $J\subseteq W_H$ containing $3$,
if $n\notin J$ or $6\notin J$, there is an independent $m$-set
$J'\subseteq W_H$ with
\[
 J'\setminus U_H=(J\setminus U_H)\setminus\{3\}.
\]
\end{lemma}
\begin{proof}
The first gap positions are $p_1=2$, $q_1=1$; hence the
chains begin $1,4,7$ and $2,5$, respectively.  The two coverage
sets both contain $3$, while neither chain does.
If $n\notin J$, replace $3$ by $2$.  The neighbors of $2$ in
$C_n^2$ are $n,1,3,4$; the independent set $J$ already avoids
$1,4$, since $1,4$ are neighbors of $3\in J$, and it avoids $n$
by hypothesis.  If $6\notin J$,
replace $3$ by $4$ instead: the neighbors of $4$ are $2,3,5,6$,
and $J$ already avoids $2,5$, since $2,5$ are neighbors of $3\in J$.  Both replacements preserve
independence and cardinality, and $2,4\in U_H\subseteq W_H$.
In either case the requirement loses exactly the free vertex $3$.
\end{proof}

\begin{theorem}[Two-gap boundary cancellation]
\label{thm:two-gap-boundary-cancellation}
For every $m\ge5$ and every admissible two-gap configuration
$H=\{h_1,h_2\}$,
\[
 h_1=3\quad\Longrightarrow\quad w_m(H)=0,
 \qquad
 h_2=2m-3\quad\Longrightarrow\quad w_m(H)=0.
\]
\end{theorem}
\begin{proof}
First suppose $h_1=3$.  By Lemma~\ref{lem:two-gap-boundary-exchange},
any requirement containing $3$ can be reduced by deleting $3$,
unless its cyclic independent set contains $\{3,6,n\}$.
We show that this exceptional triple cannot occur in an
independent $m$-set contained in $W_H$.

Write $h_2=b$.  Since $p_2>p_1=2$ and $q_2>q_1=1$, we have $b\ge5$.
If $b>5$, then $q_2\ge3$ (if $q_2=2$ then $p_2\le q_2+1=3$, so
$p_2=3$ and $b=5$), so the right chain
begins $2,5,9$.  Its coverage intervals at these vertices are
$[0,2]$, $[3,5]$, and $[7,9]$; the later coverage intervals
start beyond $9$.  Thus $6\notin W_H$, excluding the triple.

It remains to consider $H=\{3,5\}$.  Here $b=5$ forces
$(p_2,q_2)=(3,2)$, since $p_2\ge3$, $q_2\ge2$ and $p_2\le q_2+1$; hence
$P=\{2,3\}$ and $Q=\{1,2\}$.  The chains begin
$L_H=(1,4,7,10,14,18,\ldots,n-1)$ and
$R_H=(2,5,8,12,16,\ldots,n-3,n)$, with the displayed
four-step tails truncated at their indicated endpoints.
The left and right coverage windows have the exact intersection
\begin{equation}\label{eq:two-gap-boundary-coverage}
 W_H\cap[9,n-3]
 =\bigcup_{j=0}^{m-5}\{10+4j,11+4j,12+4j\}.
\end{equation}
Indeed, at each $j$ the left-chain vertex $10+4j$ contributes
$[10+4j,12+4j]$, while the right-chain vertex $12+4j$
contributes precisely the same interval, so each displayed triple
lies in $W_H$.  The vertices of $[9,n-3]$ are exactly $9$, the
displayed triples, and the vertices $13+4j$ for $0\le j\le m-6$,
since the last displayed block ends at $12+4(m-5)=n-3$.
The left intervals meeting $[9,n-3]$ are $[7,9]$ and the
$[10+4j,12+4j]$, and the right intervals meeting $[9,n-3]$ are
$[6,8]$ and the $[10+4j,12+4j]$; hence $9$ is absent from the
right coverage and each $13+4j$ belongs to neither window.
This proves \eqref{eq:two-gap-boundary-coverage}, including its
endpoints.

If an independent $m$-set $J\subseteq W_H$ contained $3,6,n$,
all its other vertices would have to lie in $[9,n-3]$.
The right side of \eqref{eq:two-gap-boundary-coverage} consists
of $m-4$ blocks of three consecutive vertices.  An independent
set of $C_n^2$ uses at most one vertex from each block.
Consequently $|J|\le3+(m-4)=m-1$, a contradiction.

We also need $\alpha(C_n^2[U_H])=m-1$.  Indeed
$a_{i+1}-b_i=A_i-\delta_i\in\{1,2\}$ for every $i<m$ (if $A_i=4$ then
$\delta_i=2$, and if $A_i=3$ then $\delta_i\in\{1,2\}$), and
$\delta_i=b_i-a_i\in\{1,2\}$; hence
$a_1,b_1,\ldots,a_m,b_m,a_1$ is a spanning $2m$-cycle in $C_n^2[U_H]$.
An independent $m$-set in this cycle is one of the alternating classes
$L_H,R_H$, neither of which is independent in $C_n^2$ since $a_1$ and
$a_m=n-1$, and $b_1=2$ and $b_m=n$, are at cyclic distance $2$.
Conversely $L_H\setminus\{a_1\}$ is independent of size $m-1$.
Thus $\alpha(C_n^2[U_H])=m-1$, so for every $T\subseteq F_H$,
$\alpha(C_n^2[U_H\cup T])=m-1$ if and only if $U_H\cup T$ contains no
independent $m$-set.

It follows that every requirement containing $3$ contains a
strictly smaller requirement not containing $3$.  In particular,
no inclusion-minimal requirement contains $3$.  Thus the
condition that $U_H\cup T$ contain an independent $m$-set is
unchanged when $3$ is toggled in $T\subseteq F_H$, and by the
preceding paragraph the same holds for the condition
$\alpha(C_n^2[U_H\cup T])=m-1$.
The fixed-point-free sign-reversing involution
$T\mapsto T\triangle\{3\}$ proves $w_m(H)=0$.

Finally, let $\varphi(i)=n+1-i$ and let $(P',Q')$ be the reflected
gap data $P'=\{m-q_2,m-q_1\}$, $Q'=\{m-p_2,m-p_1\}$.
Then $P'\subseteq\{2,\ldots,m-1\}$, $Q'\subseteq\{1,\ldots,m-2\}$, and
$q'_j\le p'_j\le q'_j+1$ for $j=1,2$ (for instance
$q'_1=m-p_2\le m-q_2=p'_1\le m-p_2+1=q'_1+1$, using
$q_2\le p_2\le q_2+1$), so $(P',Q')$ is admissible.
The gaps of $L_H$ read in reverse are the gaps of $R_{H'}$: $L_H$ has
three-gaps at $1,p_1,p_2$, so $\varphi(L_H)$, which starts at
$\varphi(a_m)=2$, has three-gaps at $m-1,m-p_1,m-p_2$, i.e.\ at
$q'_1,q'_2,m-1$.  Likewise $\varphi(R_H)=L_{H'}$.
Hence $\varphi(U_H)=U_{H'}$ and, since $\varphi$ sends a left interval
$[a,a+2]$ to the right interval $[\varphi(a)-2,\varphi(a)]$ and a right
interval $[b-2,b]$ to the left interval $[\varphi(b),\varphi(b)+2]$,
$\varphi(W_H)=W_{H'}$ and $\varphi(F_H)=F_{H'}$.
Moreover $h'_j=p'_j+q'_j=2m-h_{3-j}$, so
$H'=\{2m-h_2,2m-h_1\}$.  Since $\varphi$ is an automorphism of $C_n^2$,
it preserves independence numbers and cardinalities, so
$w_m(H)=w_m(H')$.  If $h_2=2m-3$ then $h'_1=2m-h_2=3$, and the
left-boundary case applied to $H'$ gives $w_m(H')=0$, hence $w_m(H)=0$.
\end{proof}

\subsection{A canonical block decomposition for the first interior two-gap pair}
\label{subsec:two-gap-first-interior}

We retain the notation of Subsection~\ref{subsec:two-gap-boundary}.
The following result treats the interior gap encoding $H=\{4,6\}$ for
all $m\ge5$ by a direct packing-and-uniqueness argument.

\begin{theorem}[First interior two-gap weight]
\label{thm:two-gap-46}
Let $m\ge5$, $n=4m-5$, and let $H=\{4,6\}$.  Then the
inclusion-minimal cyclic requirements in the greedy-chain interval are
exactly
\[
 \mathcal M_H=\bigl\{\{5,8,11\}\bigr\}
 \cup\bigl\{\{4j-1\}:4\le j\le m-2\bigr\}.
\]
They are pairwise disjoint, partition $F_H$, and number $m-4$.
Consequently,
\[
 \boxed{w_m(4,6)=-1.}
\]
\end{theorem}
\begin{proof}
The exceptional three-gap positions are
$P=\{2,3\}$ and $Q=\{2,3\}$.  The greedy chains are
\[
 L=(1,4,7,10,14,18,\ldots,n-1),\qquad
 R=(2,6,9,12,16,20,\ldots,n-3,n),
\]
with each four-step tail truncated at the stated endpoint.  Direct
intersection of the left and right greedy coverage windows gives
\begin{equation}\label{eq:46-coverage}
 W_H=P_0\cup E\cup D\cup\{n\},\qquad
 P_0=\{1,2,4,5,6,7,8,9,10,11,12\},
\end{equation}
where
\[
 E=\{14,16,18,\ldots,n-1\},\qquad
 D=\{15,19,\ldots,n-4\}.
\]
For $m=5$, the progression $D$ is empty.  In particular,
\[
 U_H=L\cup R,\quad |U_H|=2m,\qquad
 F_H=\{5,8,11\}\sqcup D,\quad |F_H|=m-2.
\]

First we construct cyclic independent $m$-sets realizing all the
claimed requirements.  The set
\[
 J_*=\{2,5,8,11\}\cup\{14,18,22,\ldots,n-1\}
\]
is cyclically independent, is contained in $W_H$, and satisfies
$J_*\setminus U_H=\{5,8,11\}$.  For each $4\le j\le m-2$,
put
\[
 J_j=\{2,6,9,12\}
 \cup\{16,20,\ldots,4j-4\}
 \cup\{4j-1\}
 \cup\{4j+2,4j+6,\ldots,n-1\}.
\]
The two progressions may be empty.  The four parts have respectively
$4$, $j-4$, $1$, and $m-j-1$ elements, totaling $m$.
Their consecutive cyclic gaps are at least three, including the gap
from $n-1$ back to $2$.  Thus $J_j$ is a cyclic independent $m$-set
contained in $W_H$, with
$J_j\setminus U_H=\{4j-1\}$.  This proves existence.

We next prove dominance.  We claim that $J_*$ is the \emph{unique}
cyclic independent $m$-set contained in $W_H\setminus D=P_0\cup E\cup\{n\}$.
The set $E$ has $2m-9$ vertices in steps of two and induces a path
of that order in $C_n^2$; hence $\alpha(C_n^2[E])=m-4$.
The set $P_0\subseteq[1,12]$ has independence number four.

Let $J\subseteq W_H\setminus D$ be cyclically independent of size $m$.
If $n\in J$, then $1,2,n-1\notin J$.  The available vertices
in $P_0$ lie in $[4,12]$, which contains at most three vertices
of a cyclic independent set.  Deleting $n-1$ from $E$ leaves an
even-order path on $2m-10$ vertices, with independence number $m-5$.
Therefore $|J|\le1+3+(m-5)=m-1$, a contradiction.

Hence $n\notin J$.  The bounds $|J\cap P_0|\le4$ and
$|J\cap E|\le m-4$ must both be equalities.  Since $E$ induces
an odd-order path, its maximum independent set is unique:
\[
 J\cap E=\{14,18,\ldots,n-1\}.
\]
The selected vertex $14$ excludes $12$, while $n-1$ excludes
$1$ through the cyclic wraparound.  Thus $J\cap P_0$ is an
independent four-set contained in $[2,11]$.  Its three successive
gaps must all be exactly three, so
$J\cap P_0=\{2,5,8,11\}$.  Therefore $J=J_*$,
proving uniqueness.

Now any cyclic independent $m$-set $J\subseteq W_H$ either meets
$D$, in which case its requirement contains one of the singleton
sets $\{4j-1\}$, or avoids $D$, in which case $J=J_*$ and its
requirement is $\{5,8,11\}$.  Together with the explicit
constructions above, this shows that the displayed requirements are
exactly the inclusion-minimal ones.  They are disjoint and partition
$F_H$; their number is $1+(m-5)=m-4$.

Finally, the signed avoidance sum factors over these blocks.  Each
nonempty block $M$ contributes
$\sum_{T\subsetneq M}(-1)^{|T|}=(-1)^{|M|+1}$.
Thus
\[
 w_m(4,6)=(-1)^{n-|U_H|+|F_H|+(m-4)}
 =(-1)^{(4m-5)-2m+(m-2)+(m-4)}=-1.
\]
\end{proof}

\subsection{Complete interior two-gap canonical-block theorem}
\label{subsec:complete-interior-two-gap}
The notation $L,R,U,W,F$ below agrees with $L_H,R_H,U_H,W_H,F_H$
from Subsection~\ref{subsec:two-gap-boundary}. We repeat the
coordinates to make the classification self-contained. The boundary
cancellation theorem treats $h_1=3$ or $h_2=2m-3$; the theorem below
assumes strictly interior positions and excludes the adjacent-odd
configuration. The earlier special case $H=\{4,6\}$ is retained
as an explicit example, not as a prerequisite of the general proof.
\subsubsection{Chains and free vertices}
Fix $m\ge5$, $n=4m-5$, and $P=\{p_1<p_2\}$, $Q=\{q_1<q_2\}$ with
$q_j\le p_j\le q_j+1$, $3<h_1<h_2<2m-3$, $h_j=p_j+q_j$.
Let $a_1=1,b_1=2$ and let the increments $a_{i+1}-a_i$ equal $3$ at $i=1,p_1,p_2$ and $4$ elsewhere; let $b_{i+1}-b_i$ equal $3$ at $i=q_1,q_2,m-1$ and $4$ elsewhere. Set
\[
L=\{a_1,\ldots,a_m\},\quad R=\{b_1,\ldots,b_m\},\quad U=L\cup R,
\quad W=\bigcup_i[a_i,a_i+2]\cap\bigcup_i[b_i-2,b_i]\cap[n],\quad F=W\setminus U.
\]
All intervals in the definition of $W$ are ordinary integer intervals.
Write $e_j=p_j-q_j\in\{0,1\}$ and $\delta_i=b_i-a_i$.
Direct summation gives
\[
a_i=4i-3-\mathbf1_{i>1}-\sum_j\mathbf1_{i>p_j},\quad
b_i=4i-2-\mathbf1_{i>m-1}-\sum_j\mathbf1_{i>q_j}.
\]
In particular $a_m=n-1,b_m=n$, $\delta_1=\delta_m=1$, and for $2\le i\le m-1$,
\[
\delta_i=2-\sum_{j=1}^2\mathbf1_{\{q_j<i\le p_j\}}.
\]
The two intervals in this last sum are disjoint, and each is either empty or the singleton $\{p_j\}$.

\begin{lemma}[Free vertices and canonical runs]
If $H$ is interior and nonexceptional (that is, not $h_1,h_2$ both odd with $h_2=h_1+2$), then
\[
F=\{a_i+1:2\le i\le m-1,\ \delta_i=2\}.
\]
Its maximal runs of vertices at successive distance $3$ number exactly $m-4$.
They are all singletons if $(e_1,e_2)=(1,1)$; one doubleton and the rest singletons if $e_1+e_2=1$; and two doubletons or one tripleton (with remaining singletons) if $(e_1,e_2)=(0,0)$.
\end{lemma}
\begin{proof}
Between consecutive chain vertices the only possible free vertices are $a_i+1$ when $\delta_i=2$, and $b_i+1$ when $\delta_i=\delta_{i+1}=1$ and $a_{i+1}-a_i=3$. The latter condition occurs only at $i=1,p_1,p_2$. At $i=1$ it forces $q_1=1$ and hence $h_1=3$; at $i=p_2$ it forces the excluded boundary $h_2=2m-3$; at $i=p_1$ it forces $p_1=q_1+1$, $q_2=p_1$, $p_2=p_1+1$, which is exactly the exceptional adjacent-odd case. Thus the displayed formula holds.

There are $m-2-e_1-e_2$ free vertices. A distance-$3$ link between successive free vertices occurs exactly at one of the short left-chain gaps $p_1,p_2$ with $e_j=0$. Both endpoints of that link are free: the only possible interference from the other short-gap index would contradict $q_1<q_2$; at the right endpoint the interior hypothesis excludes $p_2=m-1$ when $e_2=0$. Thus there are $2-e_1-e_2$ links, so the number of runs is $(m-2-e_1-e_2)-(2-e_1-e_2)=m-4$. With only two links the listed run shapes follow.
\end{proof}

\begin{lemma}[Base packing]
$\alpha(C_n^2[U])=m-1$.
\end{lemma}
\begin{proof}
The ordered sequence $a_1,b_1,\ldots,a_m,b_m,a_1$ is a spanning $2m$-cycle in $C_n^2[U]$. An independent $m$-set in this cycle is one of the alternating classes $L,R$, neither of which is independent in $C_n^2$ because its two endpoints are cyclically at distance $2$. Conversely $L\setminus\{a_1\}$ is independent of size $m-1$.
\end{proof}

\begin{lemma}[Realization]
For each maximal spacing-three run $B=\{a_i+1,\ldots,a_j+1\}$ of $F$, where the ellipsis means consecutive chain indices, the set
\[
J_B=\{b_1,\ldots,b_{i-1}\}\cup B\cup\{a_{j+1},\ldots,a_m\}
\]
is an independent $m$-set with $J_B\setminus U=B$.
\end{lemma}
\begin{proof}
Its size is $(i-1)+(j-i+1)+(m-j)=m$. Within each of the three portions, consecutive cyclic gaps are at least $3$. At the left junction the gap is
\[
(a_i+1)-b_{i-1}=(a_i-a_{i-1})-\delta_{i-1}+1.
\]
If the first difference is $4$ this is at least $3$; if it is $3$ and $i=2$, $\delta_1=1$; if it is $3$ and $i>2$, $\delta_{i-1}=2$ would make $a_{i-1}+1$ a free vertex joined to $B$ by a $3$-step, contradicting maximality. Thus the gap is at least $3$.
At the right junction the gap is $(a_{j+1}-a_j)-1$. If $a_{j+1}-a_j=3$, then $j\in P$ and $\delta_j=2$, hence $e_\ell=0$ at that index. The next index $j+1$ is interior and has $\delta_{j+1}=2$: for $j=p_1$, the only contrary possibility $p_2=p_1+1,e_2=1$ would give $q_2=q_1$, while for $j=p_2$ the interior condition ensures $j+1\le m-1$. Then $a_{j+1}+1$ extends $B$, a contradiction. Hence this junction gap is $3$. The closing gap $n-a_m+b_1=3$ completes the check.
\end{proof}

\subsubsection{Deletion bounds}
Write $N[x]$ for the closed neighborhood in $C_n^2$.
\begin{lemma}[Singleton deletion]
If $a_{i+1}-a_i=3$ and $\delta_i=\delta_{i+1}=2$, then for $u=a_i+1$, $v=a_{i+1}+1$,
\[
\alpha(C_n^2[U\setminus N[u]])\le m-2,\quad
\alpha(C_n^2[U\setminus N[v]])\le m-2.
\]
\end{lemma}
\begin{proof}
For $u$ use the $m-2$ cliques
\[
\{b_m,a_1,b_1\},\quad
\{a_k,b_k\}\ (2\le k<i),\quad
\{b_k,a_{k+1}\}\ (i<k<m).
\]
They cover $U$ except $a_i,b_i,a_{i+1}$, all within distance $2$ of $u$.
For $v$ use
\[
\{a_m,b_m,a_1\},\quad
\{b_k,a_{k+1}\}\ (1\le k<i),\quad
\{a_k,b_k\}\ (i+2\le k<m).
\]
They cover $U$ except $b_i,a_{i+1},b_{i+1}$, all within distance $2$ of $v$. Each listed set is a clique. Every vertex in a nonsingleton run is an endpoint of a $3$-link, so the assertion applies to it.
\end{proof}

\begin{lemma}[Path-greedy principle]
On a linear interval, earliest-available greedy selection gives a maximum independent set of the induced path square. If the available vertices are consecutive portions of interlacing chains $a_1<b_1<\cdots<a_m<b_m$ with $a_{i+1}-b_i\le2$ and $b_{i+1}-b_i\ge3$ for all $i$, after first selecting $a_i$ (respectively $b_i$), greedy follows the $a$-chain (respectively $b$-chain) until the interval ends.
\end{lemma}
\begin{proof}
Replacing the first vertex of an optimal packing by the earliest available vertex cannot invalidate subsequent selections. Repeat inductively. For chain-following, $b_i-a_i\le2$ and $a_{i+1}-a_i\ge3$; likewise $a_{i+1}-b_i\le2$ and $b_{i+1}-b_i\ge3$ by hypothesis.
\end{proof}

\begin{lemma}[Paired deletion for $(0,0)$]
Suppose $P=Q=\{r,s\}$ with $2\le r<s\le m-2$. Let $u,v$ be two free vertices left after deleting one representative from each canonical run. Then
\[
\alpha(C_n^2[U\setminus(N[u]\cup N[v])])=m-3.
\]
\end{lemma}
\begin{proof}
For $(0,0)$ all interior offsets equal $2$. The nonsingleton runs are the two doubletons at indices $(r,r+1)$ and $(s,s+1)$ when $s\ge r+2$, or the tripleton at $(r,r+1,r+2)$ when $s=r+1$. Write $u=a_i+1$, $v=a_j+1$ with $i<j$. Deleting the neighborhoods leaves the two disjoint unwrapped arcs
\[
I_1=[a_i+4,a_j-2],\qquad I_2=[a_j+4,a_i+n-2].
\]
If an arc is empty it contributes zero. Since the intervening deleted positions separate the arcs by at least three cyclic positions, there are no edges between the two surviving arcs.
There are also no wrap-around edges inside either arc. For $x<y$ in
one arc with unwrapped difference $d$, the cyclic distance is
$\min(d,n-d)$. For $I_1$, $d\le a_j-a_i-6\le n-15$ (since $a_j\le n-5$,
as $j\le m-1$ and $A_{m-1}=4$, and $a_i\ge a_2=4$), so $n-d\ge15$;
for $I_2$, $n-d=(a_j-a_i+6)+(x-(a_j+4))+((a_i+n-2)-y)\ge a_j-a_i+6\ge9$.
Hence $d\ge3$ implies cyclic distance $\ge3$, so each arc induces a
path square. The chains here satisfy
Lemma~\ref{lem:five-transitions}, so $d_i\le2$ and $B_i\ge3$;
the restricted path-greedy principle therefore applies.

We verify the greedy sequences, including their endpoints, using the path-greedy principle. When a surviving arc begins at $a_k+4$, its first chain vertex is $b_{k+1}=a_k+5$ if $a_{k+1}-a_k=3$, or $a_{k+1}=a_k+4$ if that gap is $4$. For the upper endpoint $a_\ell-2$, a right-chain run ends at $b_{\ell-1}$ if $a_\ell-a_{\ell-1}=4$; a left-chain run terminates at $a_{\ell-1}$ when that vertex is the last eligible selection, regardless of whether the final left-chain gap is $3$ or $4$. In the gap-$4$ case, the right-chain endpoint $b_{\ell-1}=a_\ell-2$ is also eligible; in the gap-$3$ case it lies beyond the upper endpoint. Across the cyclic boundary the only transitions used are $b_m=n$ to $a_2+n=n+4$, and $a_m=n-1$ to $b_1+n=n+2$; the intervening vertices are within distance $2$ of the preceding choice.

For $s\ge r+2$, the following table gives the exact greedy sequences on the two arcs. The first/last indices in each range are determined by the preceding endpoint rules; all omitted indices correspond to vertices before the lower endpoint, beyond the upper endpoint, or within distance $2$ of the previous choice.
\[
\begin{array}{c|l|l|c}
(i,j)&I_1&I_2&(g_1,g_2)\\\hline
(r,s)&b_{r+1},\ldots,b_{s-1}&b_{s+1},\ldots,b_m;\ a_2+n,\ldots,a_{r-1}+n&(s-r-1,m-s+r-2)\\
(r,s+1)&b_{r+1},\ldots,b_{s-1}&a_{s+2},\ldots,a_m;\ b_1+n,\ldots,b_{r-1}+n&(s-r-1,m-s+r-2)\\
(r+1,s)&a_{r+2},\ldots,a_{s-1}&b_{s+1},\ldots,b_m;\ a_2+n,\ldots,a_r+n&(s-r-2,m-s+r-1)\\
(r+1,s+1)&a_{r+2},\ldots,a_s&a_{s+2},\ldots,a_m;\ b_1+n,\ldots,b_{r-1}+n&(s-r-1,m-s+r-2)
\end{array}
\]
The first-arc endpoints can be checked directly from
$A_r=A_s=3$, $A_k=4$ for $k\notin\{1,r,s\}$, and
$b_k=a_k+2$ for $2\le k\le m-1$.
If $i=r$, the lower endpoint $a_r+4$ excludes $a_{r+1}=a_r+3$
and admits $b_{r+1}=a_r+5$; if $i=r+1$, the lower endpoint
$a_{r+1}+4$ admits $a_{r+2}$. At the upper endpoint $a_j-2$,
$j=s$ admits $b_{s-1}=a_s-2$ (unless the first arc is empty),
while $j=s+1$ admits $a_s=a_{s+1}-3$ and excludes
$b_s=a_{s+1}-1$. The potentially empty first arc in the third
row occurs when $s=r+2$; its listed index range is then empty. The final wrapped left-chain vertex is $a_{r-1}+n$ for $i=r$ and $a_r+n$ for $i=r+1$, since $a_i-a_{i-1}\ge3$; the next left-chain vertex lies beyond the upper endpoint $a_i+n-2$. For wrapped right-chain runs, the last vertex is $b_{r-1}+n$; the inequality $b_{r-1}\le a_r-2$ follows from the four-gap before $r$ (or $b_1=a_2-2$ when $r=2$), and $b_r=a_r+2$ lies beyond the upper endpoint. For the second-arc lower endpoint, $j=s$ gives
$a_s+4$, so the first available vertex is $b_{s+1}=a_s+5$;
$j=s+1$ gives $a_{s+1}+4$, so the first available vertex is
$a_{s+2}$. The two boundary transitions are checked by their
actual distances: from $b_m=n$ the vertices $a_1+n=n+1$ and
$b_1+n=n+2$ are forbidden and $a_2+n=n+4$ is eligible;
from $a_m=n-1$ the vertex $a_1+n=n+1$ is forbidden and
$b_1+n=n+2$ is eligible. These observations verify the endpoints
in all four rows. Each total is $m-3$.

For $s=r+1$, the three free vertices are $a_r+1,a_{r+1}+1,a_{r+2}+1$. The first arc contains no chain vertex for any pair. The second-arc greedy sequences are
\[
\begin{array}{c|l|c}
(i,j)&I_2&g_2\\\hline
(r,r+1)&b_{r+2},\ldots,b_m;\ a_2+n,\ldots,a_{r-1}+n&m-3\\
(r,r+2)&a_{r+3},\ldots,a_m;\ b_1+n,\ldots,b_{r-1}+n&m-3\\
(r+1,r+2)&a_{r+3},\ldots,a_m;\ b_1+n,\ldots,b_{r-1}+n&m-3
\end{array}
\]
Indeed $a_{r+2}=a_r+6$, $b_{r+2}=a_r+8$, and $a_{r+3}=a_r+10$. The lower endpoints of the second arcs are respectively $a_r+7,a_r+10,a_r+10$, giving the listed first vertices. Their upper endpoints are $a_r+n-2,a_r+n-2,a_{r+1}+n-2=a_r+n+1$. The last wrapped left-chain vertex in the first row is $a_{r-1}+n$; the last wrapped right-chain vertex in the other two rows is $b_{r-1}+n$, since $b_r=a_r+2$ lies beyond both corresponding endpoints. In row 1 the endpoint $a_r+n-2$ equals $b_{r-1}+n$, but it is skipped: for $r>2$ it is at distance $2$ from the preceding choice $a_{r-1}+n$, and for $r=2$ it is at distance $2$ from the preceding choice $b_m=n$. In row 3 the vertex $a_r+n$ lies inside the arc but is likewise skipped, at distance $2$ from $b_{r-1}+n$; the next vertex $b_r+n=a_r+n+2$ lies beyond the endpoint $a_r+n+1$. This includes $r=2$, when the wrapped left-chain range is empty. For $r=2$ the wrapped left-chain range in the first row is
empty; for $r=m-3$ the unwrapped left-chain ranges in the last
two rows have exactly one term. All the stated bounds remain valid.
The path-greedy principle completes the proof.
\end{proof}

\subsubsection{Canonical blocks and signed weight}
\begin{theorem}[Nonexceptional interior two-gap classification]\label{thm:complete-interior-two-gap}
For every nonexceptional interior two-gap configuration, the inclusion-minimal requirements $J\setminus U$, where $J\subseteq W$ is an independent $m$-set of $C_n^2$, are exactly the $m-4$ maximal spacing-three runs $B_1,\ldots,B_{m-4}$ of $F$. Consequently
\[
w_m(H)=(-1)^{h_1+h_2+1}.
\]
\end{theorem}
\begin{proof}
Realization supplies an independent $m$-set with requirement exactly each $B_i$. Let $X$ be any transversal, choosing one vertex from each $B_i$. If $(e_1,e_2)=(1,1)$, every block is a singleton and $W\setminus X=U$. If $e_1+e_2=1$, exactly one doubleton exists, and $W\setminus X$ adds to $U$ just one vertex from that doubleton; the singleton-deletion lemma and the base-packing lemma show its independence number is $m-1$. If $(e_1,e_2)=(0,0)$, the remaining free vertices form a pair $u,v$ from the two doubletons or from the tripleton. The singleton-deletion lemma rules out independent $m$-sets containing only one of $u,v$. The paired-deletion lemma rules out an independent $m$-set containing both, because its other $m-2$ vertices would have to lie in a graph with independence number $m-3$. Hence in all cases $\alpha(C_n^2[W\setminus X])=m-1$.

Every independent $m$-set $J\subseteq W$ therefore contains a complete block: otherwise choose from each block a vertex outside $J$, obtaining a transversal disjoint from $J$, a contradiction. Since each block is realized and the blocks are disjoint, the minimal requirements are exactly the blocks.

The signed sum over subsets $T\subseteq F$ avoiding all complete blocks factors as
\[
w_m(H)=(-1)^{n-|U|}\prod_{i=1}^{m-4}\sum_{T_i\subsetneq B_i}(-1)^{|T_i|}
=(-1)^{n-|U|+|F|+m-4}.
\]
Here $|U|=2m$, $|F|=m-2-e_1-e_2$, and $n=4m-5$. Since $h_1+h_2\equiv e_1+e_2\pmod2$, this equals $(-1)^{h_1+h_2+1}$.
\end{proof}

\begin{theorem}[Exceptional adjacent-odd two-gap cancellation]
\label{thm:exceptional-two-gap-complete}
Let $m\ge6$, $n=4m-5$, and consider an interior two-gap
configuration with $H=\{2r-1,2r+1\}$, where $3\le r\le m-3$.
Then its signed Boolean-interval weight is zero:
\[
 w_m(H)=0.
\]
\end{theorem}
\begin{proof}
In the chain notation above the exceptional positions are
$P=\{r,r+1\}$ and $Q=\{r-1,r\}$.  Their offset sequence has
$\delta_r=\delta_{r+1}=1$, while $\delta_i=2$ for
$2\le i\le m-1$ outside $\{r,r+1\}$, and
$\delta_1=\delta_m=1$.
The free-vertex classification, including the exceptional
second-type vertex, gives
\[
 F=F_0\mathbin{\dot\cup}\{x\},\qquad
 F_0=\{a_i+1:2\le i\le m-1,\ i\notin\{r,r+1\}\},
 \qquad x=b_r+1=4r-2.
\]
For each $v=a_i+1\in F_0$, the set
\[
 J_v=\{b_1,\ldots,b_{i-1}\}\cup\{v\}
       \cup\{a_{i+1},\ldots,a_m\}
\]
is an independent $m$-set in $W$ with $J_v\setminus U=\{v\}$.
Indeed, all gaps within its right-chain prefix and left-chain
suffix are at least three.  Its right junction has gap
$a_{i+1}-v=A_i-1=3$: a short left-chain gap at an interior
index can occur only at $r$ or $r+1$, where $\delta_i=1$.
Its left junction has gap
$v-b_{i-1}=A_{i-1}-\delta_{i-1}+1\ge3$:
if $A_{i-1}=3$, then $i-1\in\{1,r,r+1\}$ and
$\delta_{i-1}=1$.  Finally the cyclic closing gap is
$n-a_m+b_1=3$.  Thus each singleton $\{v\}$ is a
realized minimal requirement.

The closed neighborhood of $x$ inside $U$ is exactly
\[
 N[x]\cap U=\{a_r,b_r,a_{r+1},b_{r+1}\}.
\]
These are four consecutive vertices of the spanning cycle
$a_1,b_1,\ldots,a_m,b_m,a_1$ in $C_n^2[U]$.
The remaining $2m-4$ vertices can be partitioned into
$m-2$ consecutive pairs along that spanning cycle, each a
clique.  Consequently
\[
 \alpha(C_n^2[U\setminus N[x]])\le m-2.
\]
If an independent $m$-set $J\subseteq W$ contains $x$,
then $J\setminus\{x\}$ contains $m-1$ vertices outside
$N[x]$.  It cannot lie entirely in $U$, by the preceding
inequality.  Therefore it contains some $v\in F_0$.
Since $\{v\}$ is already a realized requirement, the
requirement $J\setminus U$ containing $x$ cannot be
inclusion-minimal.

It follows that every inclusion-minimal requirement is a
subset of $F_0$, and hence none contains $x$.  Equivalently,
for every $T\subseteq F_0$, the condition that $U\cup T$
contains an independent $m$-set is unchanged by adjoining
$x$: any such set using $x$ also uses some $v\in F_0$,
for which $J_v\subseteq U\cup T$ is already available.
Thus the deficit indicator is invariant under the toggle
$T\leftrightarrow T\cup\{x\}$.  The two terms have
opposite signs, and the signed weight vanishes.
\end{proof}

\begin{remark}
The exceptional cancellation is a direct instance of the
inessential-free-vertex mechanism.  Unlike the general
five-transition classification, the proof above uses only
the explicit two-gap coordinates and a clique cover of the
forced chain vertices.  No computational verification is
needed for the theorem.
\end{remark}

\begin{lemma}[Two-gap decomposition for $n=4m-5$]
\label{lem:two-gap-decomposition-4m-5}
For every integer $m\ge5$, put $n=4m-5$. Then
\[
 [t^{m-1}]B_{4m-5}^{\mathrm{both}}(t)=\sum_{H} w_m(H),
\]
the sum over all admissible two-gap data $(P,Q)$ as in
Subsection~\ref{subsec:two-gap-boundary}, where $w_m(H)$ is defined
in that subsection.
\end{lemma}
\begin{proof}
By the definition of $B_n^{\mathrm{both}}(t)$ and the $X=\ell-1$ argument
in Corollary~\ref{cor:both-active-3k-plus-two-top}, the coefficient is
$\sum_S(-1)^{n-|S|}$ over $S$ with $U(S)=V(S)=m>X(S)$ and $X(S)=m-1$,
and $X(S)=\alpha(C_n^2[S])$: the patterns $(1,n-1)$ and $(2,n)$ violate
the wraparound edges $\{1,n-1\}$, $\{2,n\}$ of $C_n^2$, the pattern
$(1,n)$ is impossible in $G_n$, and every other boundary pattern is
$C_n^2$-independent.
Let $S$ contribute, with greedy chains $L=(a_1<\cdots<a_m)$ and
$R=(b_1<\cdots<b_m)$ of $P_n^2[S]$. Endpoint rigidity
(Corollary~\ref{cor:both-active-endpoint-rigidity}) gives
$a_1=1$, $a_2=4$, $b_1=2$, $b_{m-1}=n-3$, $a_m=n-1$, $b_m=n$.
The three-or-four-gap restriction
(Theorem~\ref{thm:five-gap-exclusion}) gives $A_j,B_j\in\{3,4\}$,
and $\sum_{j=1}^{m-1}(A_j-3)=n-3m+1=m-4$, so each chain has exactly
three three-gaps. The first $L$-gap and the last $R$-gap are $3$,
so $L$ has exactly two further three-gaps $p_1<p_2$ and $R$ exactly
two further three-gaps $q_1<q_2$ (excluding $m-1$).
We now prove the interlacing inequalities
\[
 q_j\le p_j\le q_j+1\qquad(j=1,2)
\]
directly, without a case distinction on transition counts.
For $2\le i\le m-1$, the left chain has a three-gap at
index $1$ and at the indices $p_1,p_2$, whereas the right
chain has three-gaps at $q_1,q_2$ and $m-1$. Summing the
gaps from the initial vertices $a_1=1$ and $b_1=2$ gives
\[
 a_i=4i-4-\#\{j:p_j<i\},\qquad
 b_i=4i-2-\#\{j:q_j<i\}.
\]
Consequently
\[
 \delta_i=b_i-a_i
 =2+\#\{j:p_j<i\}-\#\{j:q_j<i\}
 \quad(2\le i\le m-1).
\]
By the greedy-offset bound, $\delta_i\in\{1,2\}$. In terms of
prefix counts $P(i)=\#\{j:p_j<i\}$ and
$Q(i)=\#\{j:q_j<i\}$, the identity reads
$Q(i)-P(i)=2-\delta_i\in\{0,1\}$.
Thus, as we scan gap positions from left to right, the right
chain has encountered either the same number of extra short
gaps as the left chain or exactly one more. This is the only
counting fact needed below.
In particular $q_j\le p_j$: otherwise at $i=p_j+1$
there would be at least $j$ positions from $P$ but at
most $j-1$ positions from $Q$, giving $\delta_i\ge3$.
If $p_j=m-1$, the inequality $q_j\le p_j$ is immediate
from $q_j\le m-2$. Otherwise $p_j\le m-2$, so
$i=p_j+1\le m-1$ lies in the range where the offset
identity applies.

Suppose instead that $p_j\ge q_j+2$. At the index
$i=q_j+1$, exactly $j$ of the $q$-positions precede $i$:
$q_j<i$ while $q_{j+1}\ge i$ if $j=1$, since $q_2>q_1$.
Exactly $j-1$ of the $p$-positions precede $i$.
For $j=1$ this follows from $i=q_1+1<p_1$.
For $j=2$, the already established inequality
$p_1\le q_1+1\le q_2<i$ shows that $p_1<i<p_2$.
Hence $\delta_i=1$.
Moreover $i\notin\{1,p_1,p_2\}$: it is greater than $1$,
lies strictly before $p_j$, and, for $j=2$, follows $p_1$.
Therefore the left-chain gap $A_i=a_{i+1}-a_i$ equals
$4$. But then
\[
 d_i=a_{i+1}-b_i=A_i-\delta_i=3,
\]
contradicting the cyclic splicing inequality
$d_i\le2$ (Lemma~\ref{lem:cyclic-splicing}).
Thus $p_j\le q_j+1$, as claimed.

Thus $e_j=p_j-q_j\in\{0,1\}$, and the contributing chain pair is
exactly one admissible two-gap datum.
Finally, the sets $S$ with given greedy chains $(L,R)$ are exactly the
Boolean interval $U_H\subseteq S\subseteq W_H$ (proof of
Theorem~\ref{thm:greedy-boundary}); distinct $(P,Q)$ give disjoint
intervals; $U(S)=V(S)=m$ throughout each interval, since $L$ and $R$
are independent $m$-sets in $G_n[S]$ with endpoint pairs $(1,n-1)$
and $(2,n)$; and the deficit condition inside the interval is exactly
$\alpha(C_n^2[U_H\cup T])=m-1$. Hence each interval contributes exactly
$w_m(H)$, and summing over intervals gives the claimed identity.
\end{proof}

\begin{corollary}[Third first-staircase coefficient]
\label{cor:staircase-4m-5}
For every integer $m\ge5$, put $n=4m-5$. Then
\[
 [t^{m-1}]B_{4m-5}^{\mathrm{both}}(t)=-1.
\]
\end{corollary}
\begin{proof}
The boundary two-gap configurations have zero weight by
Theorem~\ref{thm:two-gap-boundary-cancellation}, and the exceptional
interior configurations have zero weight by
Theorem~\ref{thm:exceptional-two-gap-complete} (for $m=5$ there is
no exceptional interior configuration, since $3<2r-1$ and
$2r+1<7$ are incompatible). For each remaining
interior configuration, Theorem~\ref{thm:complete-interior-two-gap}
gives $w_m(H)=(-1)^{h_1+h_2+1}$.

Write $p_1<p_2$, $q_i=p_i-e_i$, $e_i\in\{0,1\}$.
The interior inequalities are $3<h_1<h_2<2m-3$, where
$h_i=2p_i-e_i$. We pair configurations by replacing
$e_1$ with $1-e_1$, keeping $p_1,p_2,e_2$ fixed.
Every admissible configuration satisfies $h_2-h_1\ge2$:
indeed $h_2-h_1=2(p_2-p_1)-(e_2-e_1)$; if $p_2-p_1\ge2$
this is at least $3$, and if $p_2=p_1+1$ it is $2-(e_2-e_1)$,
which equals $1$ only when $e_2=1$, $e_1=0$, in which case
$q_1=p_1=q_2$, contradicting admissibility.
For $p_1\ge3$ the operation preserves $q_1<q_2$: if
$p_2\ge p_1+2$ then $q_1\le p_1<p_2-1\le q_2$; if $p_2=p_1+1$
and $e_2=0$ then $q_2=p_1+1>q_1$; if $p_2=p_1+1$ and $e_2=1$
then $q_2=p_1$ and admissibility holds exactly when $e_1=1$,
which is the exceptional adjacent-odd configuration
$h_1=2p_1-1$, $h_2=2p_1+1$. In particular the $e_1$-flip of an
exceptional configuration is inadmissible, so no nonexceptional
configuration pairs into an exceptional one. The partner of a
nonexceptional interior configuration with $p_1\ge3$ is again
nonexceptional and interior: $h_1'=2p_1-(1-e_1)\ge5$,
$h_1'\le h_1+1\le h_2-1$ by $h_2-h_1\ge2$, and the partner is not
exceptional, since for $e_1=0$ exceptionality would require
$h_2-h_1=1$, while for $e_1=1$ the partner has $h_1'$ even.
Hence Theorem~\ref{thm:complete-interior-two-gap} applies to both
members, and changing $e_1$ changes the parity of $h_1+h_2$ and
therefore reverses the signed weight. The only remaining failure
of interiority is $p_1=2,e_1=0$, whose partner has $h_1=3$ and is
a boundary configuration of weight zero.

For these unpaired configurations $h_1=4$, and $h_2$
runs once through the integers $6,7,\ldots,2m-4$:
the even values correspond to $e_2=0$, and the odd
values to $e_2=1$. Indeed $h_2=5$ is impossible, since
$2p_2-e_2=5$ forces $(p_2,e_2)=(3,1)$ and hence
$q_2=2=q_1$; and for $6\le h_2\le2m-4$ the unique
$(p_2,e_2)$ with $2p_2-e_2=h_2$ is $p_2=\lceil h_2/2\rceil$,
$e_2=2p_2-h_2$, which satisfies $q_2=p_2-e_2\ge3>2=q_1$
and $p_2\le m-1$. Their total signed contribution is
\[
 \sum_{h_2=6}^{2m-4}(-1)^{h_2+1}
 =-(m-4)+(m-5)=-1.
\]
By Lemma~\ref{lem:two-gap-decomposition-4m-5}, the sum over all
two-gap intervals is the coefficient, and the pairing above exhausts
all of them except the residual family just evaluated. This proves
the formula.
\end{proof}

\subsection{Cyclic deficit and Boolean finite differences for interlacing chains}
\label{subsec:universal-cyclic-deficit}

The following argument applies uniformly to interlacing chains with
any admissible number of additional short gaps.  It does not require a
classification of the minimal cyclic requirements.

\begin{theorem}[Cyclic deficit of interlacing greedy chains]
\label{thm:interlacing-cyclic-deficit}
Let $n\ge5$, and let
$L=(a_1,\ldots,a_m)$ and $R=(b_1,\ldots,b_m)$ be
independent $m$-sets of the path square $P_n^2$ satisfying
\[
 a_1=1,\qquad a_2=4,\qquad a_m=n-1,\qquad
 b_1=2,\qquad b_m=n,
\]
and, for all applicable indices,
\[
 1\le b_i-a_i\le2,\qquad
 1\le a_{i+1}-b_i\le2.
\]
Put $U=L\cup R$.  Then
\[
 \alpha(C_n^2[U])=m-1.
\]
\end{theorem}
\begin{proof}
The interlacing inequalities give
$a_1<b_1<a_2<b_2<\cdots<a_m<b_m$, with every
successive difference equal to $1$ or $2$.  As $b_m=n$
and $a_1=1$, the induced cycle-square graph on $U$
contains the spanning even cycle
\[
 a_1,b_1,a_2,b_2,\ldots,a_m,b_m,a_1.
\]
An independent set of this $2m$-cycle has size at most $m$,
and its only independent sets of size $m$ are its two
alternating classes $L$ and $R$.  Neither is independent
in $C_n^2$: the pair $\{1,n-1\}\subseteq L$ has
cyclic distance $2$, as does $\{2,n\}\subseteq R$.
Consequently $\alpha(C_n^2[U])\le m-1$.

On the other hand, $L\setminus\{1\}$ is independent in
$C_n^2$.  Its successive path gaps are at least $3$,
while its wraparound gap from $a_m=n-1$ to $a_2=4$
is $n-(n-1)+4=5$.  It has $m-1$ vertices, giving
the reverse inequality.
\end{proof}

\begin{corollary}[Two-value property on the greedy-chain interval]
\label{cor:interlacing-two-value}
Under the hypotheses of Theorem~\ref{thm:interlacing-cyclic-deficit},
suppose $L,R$ are the left- and right-greedy maximum
independent chains of their associated Boolean interval.  Set
\[
 W=\left(\bigcup_{i=1}^m[a_i,a_i+2]\right)
 \cap\left(\bigcup_{i=1}^m[b_i-2,b_i]\right)\cap[n],
 \qquad F=W\setminus U.
\]
Then for every $T\subseteq F$,
\[
 \alpha(C_n^2[U\cup T])\in\{m-1,m\}.
\]
The assertion means that at most these two values occur;
it does not assert that both are attained.
\end{corollary}
\begin{proof}
The greedy-chain Boolean-interval property gives
$\alpha(P_n^2[U\cup T])=m$, so the cycle-square
independence number is at most $m$.  Since $U\subseteq U\cup T$,
Theorem~\ref{thm:interlacing-cyclic-deficit} gives the
lower bound $m-1$.
\end{proof}

For a set function $f$ on subsets of $[n]$, write
\[
 \Delta_F f(U)=\sum_{T\subseteq F}
 (-1)^{|F|-|T|}f(U\cup T).
\]

\begin{corollary}[Signed weight as a Boolean finite difference]
\label{cor:interlacing-finite-difference}
Under the hypotheses of Corollary~\ref{cor:interlacing-two-value},
let $\alpha_{\mathrm{cyc}}(S)=\alpha(C_n^2[S])$ and
\[
 w(L,R)=\sum_{\substack{T\subseteq F\\
 \alpha_{\mathrm{cyc}}(U\cup T)=m-1}}
 (-1)^{n-|U|-|T|}.
\]
If $F\ne\varnothing$, then
\[
 \boxed{\quad
 w(L,R)=(-1)^{n-|U|+|F|+1}
 \Delta_F\alpha_{\mathrm{cyc}}(U)
 =(-1)^{n+|F|+1}\Delta_F\alpha_{\mathrm{cyc}}(U).
 \quad}
\]
If $F=\varnothing$, then $w(L,R)=(-1)^{n-2m}$.
\end{corollary}
\begin{proof}
By the two-value property, the indicator of the event
$\alpha_{\mathrm{cyc}}(U\cup T)=m-1$ is
$m-\alpha_{\mathrm{cyc}}(U\cup T)$.  Thus
\[
 w(L,R)=(-1)^{n-|U|}
 \sum_{T\subseteq F}(-1)^{|T|}
 \bigl(m-\alpha_{\mathrm{cyc}}(U\cup T)\bigr).
\]
When $F$ is nonempty, the constant term vanishes because
$\sum_{T\subseteq F}(-1)^{|T|}=0$; rearranging signs
proves the first formula.  The interlacing inequalities
make the chains disjoint, hence $|U|=2m$, giving the
second formula.  If $F$ is empty, $U$ alone contributes,
and its cyclic independence number is $m-1$ by the theorem.
\end{proof}

\section{Forward-slots decompositions and M\"obius distributions for cycle powers}
\label{sec:third-arrangement}
For integers $r\ge1$ and $n\ge2r+1$, write $C_n^r$ for the graph on
$\mathbb Z_n$ in which two distinct vertices are adjacent when their
cyclic distance is at most $r$. For an independent set $I$ of $C_n^r$,
define
\[
 Q_r(I)=\bigcup_{x\in I}\{x,x+1,\ldots,x+r\}\subseteq\mathbb Z_n.
\]
The case $r=2$ recovers the cycle-square forward slots used above.

\begin{theorem}[Forward-slots bijection for cycle powers]
\label{thm:slots-bijection}
Let $r,k\ge1$ and $n\ge2rk+r+1$. Every $S\subseteq\mathbb Z_n$
with $\alpha(C_n^r[S])=k$ admits exactly one independent $k$-set
$I\subseteq S$ such that $S\subseteq Q_r(I)$. Consequently, for
$k\le j\le(r+1)k$,
\[
 \#\{S\subseteq\mathbb Z_n:|S|=j,\ \alpha(C_n^r[S])=k\}
 =i_k^{(r)}(n)\binom{rk}{j-k},\qquad
 i_k^{(r)}(n)=\frac{n}{n-rk}\binom{n-rk}{k}.
\]
\end{theorem}
\begin{proof}
Suppose first that all clockwise gaps between consecutive elements
of $S$ are at most $r$. Starting at any $x\in S$, repeatedly choose
the first vertex of $S$ at clockwise distance at least $r+1$ from
the previous choice. Each step has length at most $2r$: otherwise
$r$ consecutive positions at distances $r+1,\ldots,2r$ would be
missing from $S$. After $k$ steps the chosen $k+1$ vertices have
successive gaps at least $r+1$ and lie in a clockwise arc of length
at most $2rk$. Since $n-2rk\ge r+1$, their closing gap is also
at least $r+1$. This gives an independent $(k+1)$-set, a
contradiction. Hence $S$ has a clockwise gap of length at least $r+1$.

Cut the cycle at such a gap and order the available vertices along
the resulting path. Any path-independent set is also cyclically
independent, because the cut gap has length at least $r+1$.
The earliest-available greedy algorithm selects a maximum independent
set $I$ of size $k$; the standard exchange argument replaces the
first vertex of any independent set by the earliest available vertex,
then proceeds inductively. Every available vertex is within $r$
clockwise positions of its preceding greedy selection, including
those after the last selection. Therefore $S\subseteq Q_r(I)$.

For uniqueness, suppose $I$ and $J$ are two representations. The $k$ forward slots of $J$ are pairwise disjoint: two slot
starts are at clockwise distance at least $r+1$, so their
length-$(r+1)$ clockwise intervals do not overlap. Each slot
induces a clique in $C_n^r$. Since the independent $k$-set
$I$ is contained in the union of these $k$ cliques, it has at
most one vertex in each slot; counting $k$ vertices forces
exactly one vertex in every slot. List the slot
starts $y_1,\ldots,y_k$ clockwise and their matched $I$-vertices
$x_1,\ldots,x_k$, in the same cyclic order: the slots are disjoint
clockwise arcs, and each contains exactly one matched vertex, so
matching cannot reverse or permute their cyclic order. Write
$d_i=x_i-y_i\in\{0,\ldots,r\}$ in the clockwise lift.
If $d_i>0$, then $y_i$ cannot belong to the $I$-slot beginning
at $x_i$; as $I$ has exactly one vertex in each $J$-slot, the
only possible $I$-slot containing $y_i$ is the preceding one.
Indeed, all $I$-slot starts are in the same cyclic order as the
matched $x_i$, and the next start follows $y_i$. Thus
$0\le y_i-x_{i-1}\le r$. If $d_{i-1}=0$, this would put
consecutive starts $y_{i-1},y_i$ of $J$ within $r$ positions,
contradicting independence. Hence $d_{i-1}>0$. Iterating
backwards shows that every $d_i>0$, and therefore
\[
 x_i-x_{i-1}=(x_i-y_i)+(y_i-x_{i-1})\le2r
\]
for every cyclic index. Summing these $k$ cyclic gaps yields
$n\le2rk$, contrary to the hypothesis. Thus all $d_i=0$,
and $I=J$.

For each independent $I$, the $k$ disjoint forward slots cover
$(r+1)k$ vertices, of which $k$ belong to $I$. Hence the
Boolean interval $[I,Q_r(I)]$ has $rk$ optional vertices.
Uniqueness makes these intervals disjoint. Finally, counting
clockwise gaps of length at least $r+1$ between successive
vertices of $I$ gives
\[
 i_k^{(r)}(n)=\frac nk\binom{n-rk-1}{k-1}
 =\frac{n}{n-rk}\binom{n-rk}{k}.
\]
\end{proof}

\subsection{Unique reconstruction of partial forward slots}
\label{subsec:partial-slots}
For $U\subseteq\mathbb Z_n$, put
\[
 E_{n,k}^{(r)}(U)=\sum_{\substack{U\subseteq S\subseteq\mathbb Z_n\\
 \alpha(C_n^r[S])=k}}(-1)^{|S|-|U|}.
\]
This is the upper Boolean-lattice M\"obius transform of the
independence-number-$k$ level indicator. We write
$E_{n,k}(U)=E_{n,k}^{(2)}(U)$ when $r=2$.

\begin{lemma}[Partial-slot reconstruction for cycle powers]
\label{lem:partial-slot-injectivity}
Let $r,k\ge1$ and $n>(r+1)k$. If $I$ is an independent
$k$-set of $C_n^r$ and $J\subseteq I$, then
$Q_r(I)\setminus J$ uniquely determines both $I$ and $J$.
\end{lemma}
\begin{proof}
The forward slots of $I$ are disjoint. After deleting $J$,
each slot contributes either a full block of $r+1$ consecutive
vertices or a shortened block of $r$ consecutive vertices
obtained by deleting its starting vertex. Put $T=Q_r(I)\setminus J$.
Since $|T|\le(r+1)k<n$, the maximal consecutive clockwise runs
of $T$ have well-defined first vertices.

If two successive blocks are consecutive without an intervening
missing vertex, the second block must be full. Indeed, a
shortened second block beginning immediately after the preceding
block would have its underlying slot start at the preceding
block's final vertex. Its start would then be at distance at most
$r$ from the preceding slot start, contrary to independence.
Thus each maximal run consists either entirely of full blocks,
or of one initial shortened block followed by full blocks.
Its length is consequently congruent to $0$ or $r$ modulo $r+1$.
These residues are distinct, and they uniquely determine the
block decomposition from the first vertex of the run. A full
block begins at its slot start; a shortened block begins one
vertex after its slot start. This reconstructs $I$ and identifies
exactly the deleted starts $J$.
\end{proof}

\begin{theorem}[Complete stable-range M\"obius distribution for cycle powers]
\label{thm:stable-mobius-distribution}
Let $r,k\ge1$ and $n\ge2rk+r+1$. For every
$U\subseteq\mathbb Z_n$ with $\alpha(C_n^r[U])=k$,
\[
 E_{n,k}^{(r)}(U)=
 \begin{cases}
 (-1)^{|J|},&U=Q_r(I)\setminus J\text{ for an independent }k\text{-set }I
                    \text{ and }J\subseteq I,\\
 0,&\text{otherwise}.
 \end{cases}
\]
The pair $(I,J)$ in the first case is unique. In particular,
for $0\le s\le k$,
\[
 \#\{U:|U|=rk+s,\ \alpha(C_n^r[U])=k,\
 E_{n,k}^{(r)}(U)=(-1)^{k-s}\}
 =\frac{n}{n-rk}\binom{n-rk}{k}\binom{k}{s}.
\]
All remaining sets of independence number $k$ have coefficient
zero. The number of nonzero coefficients is
$2^k i_k^{(r)}(n)$, with $2^{k-1}i_k^{(r)}(n)$ positive
and the same number negative.
\end{theorem}
\begin{proof}
By Theorem~\ref{thm:slots-bijection}, the level-$k$ sets form
the disjoint Boolean intervals $[I,Q_r(I)]$. Fix one such
interval and write $Q=Q_r(I)$. Its contribution is empty if
$U\nsubseteq Q$. Otherwise every summand has the form
$S=I\cup U\cup T$, where $T$ ranges over the subsets of
$Q\setminus(I\cup U)$. Therefore its signed contribution is
\[
 (-1)^{|I\setminus U|}
 \sum_{T\subseteq Q\setminus(I\cup U)}(-1)^{|T|}.
\]
The inner sum is $(1-1)^{|Q\setminus(I\cup U)|}$: it vanishes
unless $Q=I\cup U$, when the contribution is
$(-1)^{|I\setminus U|}$. This condition is equivalent to
$U=Q_r(I)\setminus J$, with $J=I\setminus U\subseteq I$.
Lemma~\ref{lem:partial-slot-injectivity} shows that at most
one pair contributes. Conversely, for every $J\subseteq I$,
the translated set $I+r$ is independent and lies in
$Q_r(I)\setminus J$. Since this latter set is covered by $k$
disjoint clique slots, its independence number is exactly $k$.
If $|U|=rk+s$, then $|J|=k-s$. There are $\binom{k}{s}$
choices of $J$ for each of the $i_k^{(r)}(n)$ choices of $I$,
and all images are distinct. Summation over $s$ and parity
of $k-s$ proves the remaining assertions.
\end{proof}

\subsection{M\"obius inversion on restricted Boolean intervals}
\label{subsec:restricted-mobius}
For $U\subseteq W\subseteq\mathbb Z_n$ and $a=\alpha(C_n^2[U])$,
write
\[
 E_a(U,W)=\sum_{\substack{U\subseteq S\subseteq W\\
               \alpha(C_n^2[S])=a}}(-1)^{|S|-|U|}.
\]
The full-upper-endpoint coefficients in Theorem~\ref{thm:stable-mobius-distribution}
are the special case $W=\mathbb Z_n$.

\begin{lemma}[Boolean-interval alternating sum]
\label{lem:interval-transform}
If $L\subseteq R$ and $U\subseteq W$ are subsets of a finite set, then
\[
\sum_{\substack{L\subseteq S\subseteq R\\ U\subseteq S\subseteq W}}
 (-1)^{|S|-|U|}
 =\begin{cases}
 (-1)^{|L\setminus U|},&L\subseteq W,\ U\subseteq R,\
                         R\cap W=L\cup U,\\
 0,&\text{otherwise}.
 \end{cases}
\]
\end{lemma}
\begin{proof}
If $L\not\subseteq W$ or $U\not\subseteq R$, the sum is empty.
Otherwise its indexing interval is $[L\cup U,R\cap W]$,
and its signed sum is
$(-1)^{|L\setminus U|}(1-1)^{|(R\cap W)\setminus(L\cup U)|}$.
\end{proof}

For this subsection, write $Q(I)=Q_2(I)$ for the forward-slot
union in the square of a cycle.

\begin{theorem}[Restricted forward-slots transform]
\label{thm:restricted-mobius}
Let $k\ge1$, $n\ge4k+3$, and $U\subseteq W\subseteq\mathbb Z_n$
with $\alpha(C_n^2[U])=k$. Then
\[
\boxed{ E_k(U,W)=
\sum_{\substack{I\in\mathcal I_k(C_n^2),\ I\subseteq W\\
 U\subseteq Q(I),\ Q(I)\cap W=U\cup I}}
 (-1)^{|I\setminus U|}. }
\]
\end{theorem}
\begin{proof}
Theorem~\ref{thm:slots-bijection} partitions the independence-level-$k$
sets into the disjoint Boolean intervals $[I,Q(I)]$.
Apply Lemma~\ref{lem:interval-transform} to each interval with
$L=I$ and $R=Q(I)$, and sum the resulting contributions.
\end{proof}

\subsection{Deficit complexes and their algebraic structure}
\label{subsec:deficit-complexes}
The following constructions apply to an arbitrary finite graph $G$.
For $U\subseteq W\subseteq V(G)$, let $a=\alpha(G[U])$ and
$F=W\setminus U$. Define the simplicial complex
\[
 \Delta_G(U,W)=\{T\subseteq F:\alpha(G[U\cup T])=a\}.
\]
The complex contains the empty face, and is downward closed since
$\alpha(G[U\cup T])\ge a$ for all $T$.
Write $E_G(U,W)=\sum_{T\in\Delta_G(U,W)}(-1)^{|T|}$.
With the augmented convention for reduced Euler characteristic,
\begin{equation}\label{eq:euler-deficit}
 E_G(U,W)=-\widetilde\chi(\Delta_G(U,W)).
\end{equation}
Here $\widetilde\chi(\Delta)=\sum_{T\in\Delta}(-1)^{|T|-1}$,
including the empty face; this convention also covers the complex
$\{\varnothing\}$.

\begin{proposition}[Factorization over disconnected induced graphs]
\label{prop:deficit-factorization}
Suppose $W=W_1\mathbin{\dot\cup}\cdots\mathbin{\dot\cup}W_s$
and no edge of $G$ joins distinct $W_i$. Put $U_i=U\cap W_i$.
Then
\[
 \Delta_G(U,W)=\Delta_G(U_1,W_1)*\cdots*
 \Delta_G(U_s,W_s),\qquad
 E_G(U,W)=\prod_{i=1}^s E_G(U_i,W_i).
\]
\end{proposition}
\begin{proof}
Every $T\subseteq W\setminus U$ decomposes uniquely into
$T_i\subseteq W_i\setminus U_i$. Independence numbers add across
disjoint induced components, and each local independence number is
at least $\alpha(G[U_i])$. Hence the global independence number
is unchanged exactly when every local independence number is
unchanged. This gives the join identity, and the signed face sum
factors because $(-1)^{|T|}=\prod_i(-1)^{|T_i|}$.
\end{proof}

\begin{proposition}[Minimal nonfaces and independent witnesses]
\label{prop:deficit-nonfaces}
Every inclusion-minimal nonface $T$ of $\Delta_G(U,W)$ has
$|T|\le a+1$. More precisely, there exists an independent
$(a+1)$-set $A\subseteq U\cup T$ such that $T=A\setminus U$.
\end{proposition}
\begin{proof}
Because $T$ is a nonface, $U\cup T$ contains an independent
$(a+1)$-set $A$. Then $A\setminus U\subseteq T$ is itself a
nonface. Minimality forces $A\setminus U=T$.
\end{proof}

\begin{corollary}[Quadratic--cubic presentation at level two]
\label{cor:level-two-ideal}
Suppose $a=2$ and every $x\in F$ is individually safe, that is,
$\alpha(G[U\cup\{x\}])=2$. Let $B$ be the graph on $F$ whose
edges are pairs $\{x,y\}$ with
$\alpha(G[U\cup\{x,y\}])>2$.
Let $\mathcal H$ be the family of independent triples of $G[F]$
containing no edge of $B$. Then the Stanley--Reisner ideal of
$\Delta_G(U,W)$ over any field $\Bbbk$ is
\[
 I_{\Delta}=\bigl(X_xX_y:\{x,y\}\in E(B)\bigr)
 +\bigl(X_xX_yX_z:\{x,y,z\}\in\mathcal H\bigr).
\]
Without the individual-safety hypothesis, the ideal is still
generated in degrees at most three, but linear generators may occur.
\end{corollary}
\begin{proof}
Proposition~\ref{prop:deficit-nonfaces} bounds the size of a
minimal nonface by three. Individual safety excludes singleton
nonfaces. The minimal nonfaces of size two are exactly the edges
of $B$. If all pairs of a triple $T$ are faces, an independent
three-set newly appearing in $U\cup T$ must contain all three
vertices of $T$: otherwise a pair would already be a nonface.
Thus $T$ is a minimal nonface precisely when it is independent
in $G[F]$ and contains no edge of $B$. The displayed generators
are exactly the squarefree monomials associated with minimal
nonfaces.
\end{proof}

\subsection{Canonical requirements and topological cancellation}
\label{subsec:requirements-topology}
The minimal-requirement method in
Proposition~\ref{prop:canonical-requirement-cover} also admits a topological
interpretation. The following statement isolates the exact
hypothesis needed for a sphere-or-cone conclusion.

\begin{proposition}[Disjoint requirements and joins]
\label{prop:requirements-join}
Let $\Delta$ be a simplicial complex on a finite vertex set $F$.
Suppose its inclusion-minimal nonfaces are exactly the pairwise
disjoint nonempty sets $B_1,\ldots,B_s\subseteq F$.
Put $P=F\setminus\bigcup_i B_i$. Then
\[
 \Delta=\left(*_{i=1}^s\partial\sigma(B_i)\right)*\sigma(P),
\]
where $\sigma(P)$ is the simplex on $P$ and
$\partial\sigma(B_i)$ is the boundary complex of the simplex
on $B_i$ (with the augmented convention for $|B_i|=1$).
In particular, if $P\ne\varnothing$, then $E(\Delta)=0$ and
$\Delta$ is contractible. If $P=\varnothing$, then
\[
 E(\Delta)=(-1)^{|F|+s},
\]
and $\Delta$ has the augmented homotopy type of a sphere of
dimension $\sum_i(|B_i|-1)-1$.
\end{proposition}
\begin{proof}
A subset of $F$ is a face exactly when it contains none of
the $B_i$ in its entirety. Because the blocks are disjoint,
face choices factor into proper subsets of each $B_i$ and
arbitrary subsets of $P$, proving the join identity. For each nonempty $B_i$, the signed sum over its proper subsets is
\[
\sum_{T\subsetneq B_i}(-1)^{|T|}
= -(-1)^{|B_i|}=(-1)^{|B_i|+1}.
\]
Thus the product of the block contributions is
$(-1)^{\sum_i(|B_i|+1)}=(-1)^{|F|-|P|+s}$.
Multiplying the block contributions by the unrestricted $P$-factor gives
\[
 E(\Delta)=(-1)^{|F|-|P|+s}(1-1)^{|P|}.
\]
Hence the total signed sum vanishes when $P\ne\varnothing$,
and equals $(-1)^{|F|+s}$ when $P=\varnothing$.
The homotopy assertions follow from the join description.
\end{proof}

\section{Concluding discussion}
The induced independence-number distribution is a natural
refinement of counting independent sets: it records the
largest independent sets inside every vertex-induced
subgraph. For powers of cycles, the forward-slots
decomposition makes this distribution explicit in a
stable range and gives a corresponding description of
the upper Boolean-lattice M\"obius transform. For cycle
squares, greedy chains provide a different route to
signed weights, including staircase and two-gap
cancellation identities.

These two methods complement the existing theory of
independence polynomials and induced-subgraph graph
polynomials. The independent-set counting factor is
classical, while the contribution of the forward-slots
approach is the representation of \emph{all} subsets
at a fixed induced independence level. The signed
cycle-square theory addresses cancellations not
visible from the ordinary independent-set count.
The restricted-transform and deficit-complex
formulations place these calculations in familiar
Boolean-lattice and simplicial-complex frameworks.

\section*{Declaration on the use of artificial intelligence}
Generative artificial intelligence tools, including ChatGPT, Claude,
Gemini, and Muse, were used during the development and revision of this
manuscript to assist with mathematical exploration, proof analysis and
verification, identification and correction of errors, literature
review, and improvement of the English exposition. AI-generated
arguments and suggestions were critically reviewed and revised by the
authors. The authors take full responsibility for the mathematical
statements, proofs, citations, and final manuscript.

\end{document}